\documentclass[10pt,reqno]{amsart}
\usepackage{macros}
\usepackage{xurl}

\title{Pair Siegel-Veech constants for branched covers of tori and Veech surfaces}
\author[Fenves]{Inbo Gottlieb Fenves}
\address{University of Wisconsin-Madison, Madison, WI, U.S.}
\email{gottliebfenv@wisc.edu}
\date{}

\begin{document}

\begin{abstract}
    \noindent For a locus $\Ec$ of branched covers of tori satisfying a monodromy assumption, we show there exists an explicit constant $c_A(\Ec)$ so that the number of pairs of saddle connections on surfaces in $\Ec$ of length at most $R$ and cross product at most $A$ grows like $c_A(\Ec)R^2$ for every Diophantine surface in $\Ec$; this constant is moreover piecewise linear in $A$. In particular, this covers the case of billiards on rectangles with a barrier. Similar pointwise results are shown for Veech surfaces. The proof uses work of Margulis-Mohammadi on the quantitative inhomogeneous Oppenheim conjecture.
\end{abstract}

\maketitle

    \setcounter{tocdepth}{1}
    \tableofcontents
    \setcounter{tocdepth}{2}

\section{Introduction}

Let $\Hc$ be a stratum of translation surfaces, with unit-area locus $\Hc^1 \subseteq \Hc$. Given a surface $\omega \in \Hc$, let $\Hol_\omega \subseteq \C$ denote its set of holonomy vectors. This set is discrete, and so for each $A \ge 0$ and $R > 0$ the quantity
\begin{equation*}
    N_A(\omega,R) = \left| \left\{ (z,w) \in \Hol_\omega \times \Hol_\omega : |z|^2+|w|^2\le R^2,\ |z \wedge w| \le A \right\} \right|
\end{equation*}
is finite. The function $N_A(\omega,R)$ counts the number of triangles in $\Hol_\omega$ with a vertex at the origin, of $R$-bounded lengths and area at most $A/2$. For $\Mc \subseteq \Hc$ an affine subvariety, $c_A(\Mc)$ is the unique constant for which the limit
\begin{equation}\label{eq:def-cA}
    c_A(\omega) = \lim_{R \to \infty} \frac{1}{R^2}N_A(\omega,R)
\end{equation}
exists and is equal to $c_A(\Mc)$, for $\mu_{\Mc}$-a.e. $\omega \in \Mc^1$, with $\mu_\Mc$ the volume on $\Mc^1$. It was established by Athreya-Fairchild-Masur \cite{AthFaiMas} that this constant exists\footnote{\cite{AthFaiMas} uses a slightly different choice of normalization for the counting function $N_A(\omega,R)$.} when $\Mc$ is a connected component of a stratum. On the other hand, in forthcoming joint work with Fairchild-Katipearachchi-Masur \cite{FaiFenKatMas}, we show that in any stratum $\Hc(\kappa)$ with positive rel, there exist uncountably many surfaces $\omega \in \Hc^1(\kappa)$ for which
\begin{equation}\label{eq:FFKMcounterexamples}
    N_A(\omega,R) \ge \frac{R^2\log R}{\log\log R}
\end{equation}
holds for all $A > 0$ and arbitrarily large $R \gg_A 1$; this implies the limit \eqref{eq:def-cA} need not exist pointwise.

The primary goal of this paper is to show pointwise convergence of \eqref{eq:def-cA} for a wide range of rank-1 subvarieties $\Ec$ subject to (necessary) Diophantine conditions, and in addition to explicitly compute the values $c_A(\Ec)$. 

We say a branched cover $p \colon \omega \to \eta$ is \textit{balanced} if the distinguished points (i.e., cone or marked points) of $\omega$ are exactly the preimage of the distinguished points on $\eta$. This notion was introduced by Gutkin \cite{Gut} (see also \cite{HubSch}, \cite{GutHubSch}).

\begin{thm}\label{thm:1convergenceforbranchedcovers}
    Let $\Ec$ be a subvariety consisting of (possibly disconnected) degree-$d$ balanced branched covers of tori, branched over $n \ge 2$ points, and of full rel dimension $\rel \Ec = n-1$. Then for every Diophantine\footnote{See Definitions \ref{defn:diophantinemarkedtori} and \ref{defn:diophantinesurf} for a definition.} surface $\omega \in \Ec^1$,
    \begin{equation}\label{eq:1convergenceforbranchedcovers}
        c_A(\omega) = \left( \frac{6}{\pi}d + \pi n(n-1)d \right) + 12d\Phi_3(dA) + (2\pi^2d^2n^2(n-1)^2 + 24d^2n(n-1))A,
    \end{equation}
    where $\Phi_3(T) = \sum_{t=1}^{\floor{T}} \frac{\phi(t)}{t}$, and $\phi$ is the Euler totient function.
\end{thm}

\begin{cor}\label{cor:cAEcaffinefunction}
    If $\Ec$ is as in Theorem \ref{thm:1convergenceforbranchedcovers}, then $A \mapsto c_A(\Ec)$ is a well-defined piecewise affine function, with jump discontinuities exactly at $d^{-1}\N$.
\end{cor}

Notably, Theorem \ref{thm:1convergenceforbranchedcovers} holds \textit{pointwise} for Diophantine surfaces, as opposed to on an inexplicit full-measure set. Moreover, \eqref{eq:1convergenceforbranchedcovers} is independent of the monodromy of the cover $\Ec \to \Hc(0)$, so long as the balanced assumption is satisfied. By combining Theorem \ref{thm:1convergenceforbranchedcovers} with the exceptional surfaces of \cite{FaiFenKatMas}, we obtain the following sharp dimension gap for surfaces with exceptional pair counting.

\begin{thm}\label{thm:2dimensiongapforcAEc}
    If $\Ec$ is as in Theorem \ref{thm:1convergenceforbranchedcovers}, then
    \begin{equation*}
        \dim_H \left\{ \omega \in \Ec : \text{$\inf_{A> 0}\limsup_{R \to \infty}\frac{1}{R^2}N_A(\omega,R) =\infty$} \right\} = \dim_H \Ec - 2,
    \end{equation*}
    where $\dim_H$ is the Hausdorff dimension. Moreover, this set is always uncountable.
\end{thm}

Omri Solan has recently announced a partial classification theorem for horocycle-invariant measures on moduli space, which is known to imply convergence of (single-variable) Siegel-Veech constants for \textit{every} translation surface $\omega$, with the limiting value depending only on the orbit-closure and area of the surface. Theorem \ref{thm:2dimensiongapforcAEc} says this cannot happen for the constants $c_A$.\footnote{The counterexamples of \cite{FaiFenKatMas} are branched covers of tori, and thus are ``essentially homogeneous;'' in particular, for these surfaces, exact quadratic asymptotics follows from Ratner's theorems as in Eskin-Masur-Schmoll \cite{EskMasSch}.}

Note that not all branched covers $\Ec \to \Hc(0)$ will satisfy the balanced condition of Theorem \ref{thm:1convergenceforbranchedcovers}, since preimages of branched points may be regular in a cover. However, the conditions are satisfied in the historically significant case of slit torus constructions, so we obtain the following immediate corollary.

\begin{cor}\label{cor:cAE411}
    Suppose $\omega \in \Hc(1,1)$ is a unit-area slit torus with Diophantine slit. Then
    \begin{equation}\label{eq:cAE411}
        c_A(\omega) = \frac{12}{\pi} + 4\pi + 24\Phi_3(2A) + (32\pi^2 + 192)A.
    \end{equation}
    In particular, \eqref{eq:cAE411} is the pair Siegel-Veech constant for the eigenform locus $\Ec_4(1,1)$.
\end{cor}

\begin{figure}[H]
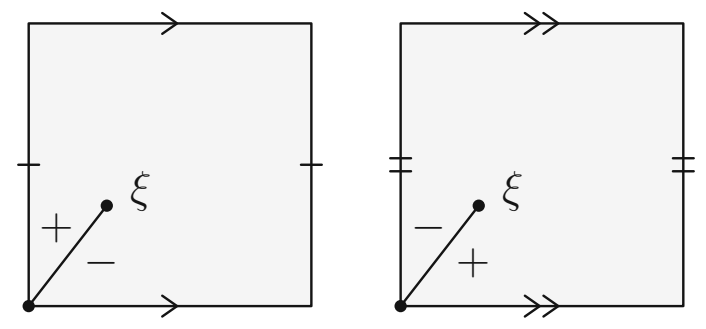

    \centering
    \SlitToriFigure[0.7]
    \caption{A slit torus construction in $\Hc(1,1)$}
    \label{fig:identicalslittori}
\end{figure}

In Section \ref{subsec:nomarkedpts}, we provide a necessary and sufficient combinatorial condition for a stratum $\Hc(\kappa)$ to contain a subvariety $\Ec$ as in Theorem \ref{thm:1convergenceforbranchedcovers}.

The primary tool we use to establish Theorems \ref{thm:1convergenceforbranchedcovers} and \ref{thm:2dimensiongapforcAEc} is the work of Margulis-Mohammadi \cite{MarMoh} on inhomogeneous versions of the quantitative Oppenheim conjecture in signature $(2,2)$. The homogeneous version was considered previously by Eskin-Margulis-Mozes \cite{EskMarMoz}, \cite{EskMarMoz1}. See also \cite{Kim} for the signature $(2,1)$-case, as well as \cite{StrSodVis} and the references therein for effective results. We remark that complementary results on two-point correlations for affine lattices were obtained by Marklof \cite{Mar} and El-Baz-Marklof-Vinogradov \cite{ElBMarVin}.

\begin{figure}[H]
    \centering
    \includegraphics[width=0.7\linewidth]{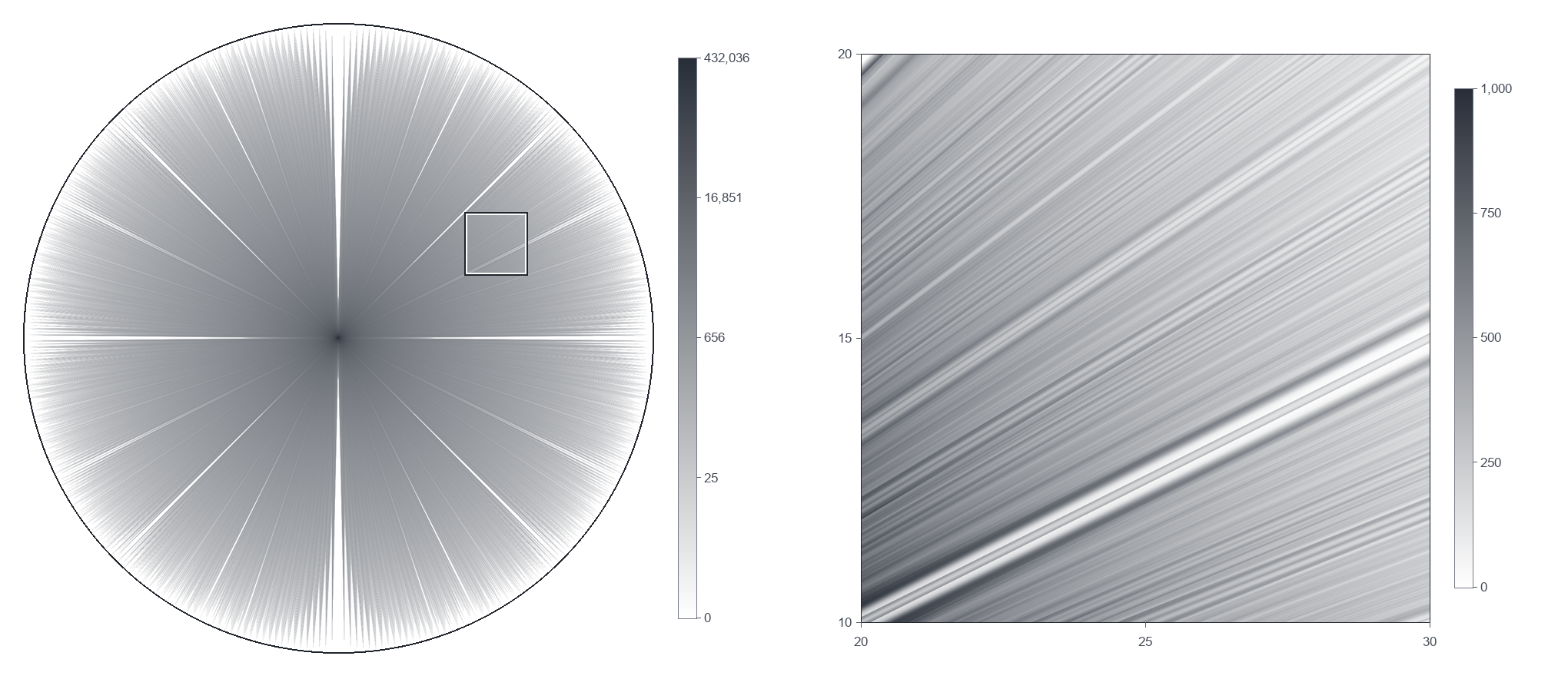}
    \caption{The distribution of triangles on $\omega \in \Ec_4(1,1)$ with $A = 10$, $R = 50$}
    \label{fig:triangledensity}
\end{figure}

As a further application of our methods, we prove pointwise convergence of $c_A(\omega)$ for every Veech surface $\omega$, and compute this constant in terms of some finite combinatorial data.

\begin{thm}\label{thm:3convergenceforveech}
    If $\omega$ is Veech, then $c_A(\omega)$ exists and depends only on its $\SL_2(\R)$-orbit. Moreover, for any Veech curve $\Vc \subseteq \Hc$, there exists a constant $\kappa_\Vc > 0$, a discrete set $\VT(\Vc) \subseteq [0,\infty)$ and positive integer-valued weights $\phi_\Vc \colon \VT(\Vc) \to \N$ for which
    \begin{equation}\label{eq:3convergenceforveech}
        c_A(\omega) = \frac{\kappa_\Vc}{\Vol(\Vc^1)} + \frac{2\pi^2}{\Vol(\Vc^1)}\sum_{t \in (0,A] \cap \VT(\Vc)} \frac{\phi_\Vc(t)}{t}
    \end{equation}
    holds for every $\omega \in \Vc^1$, where $\Vol(\Vc^1) = \Covol(\Stab_{\SL_2(\R)}(\omega))$ is the natural volume on $\Vc^1$.
\end{thm}

\begin{rmk}
    Theorem \ref{thm:3convergenceforveech} was proven in \cite[Proposition 2.1]{AthFaiMas} in the special case $A=0$, for a different choice of normalization on $\R^4$; the more general Theorem \ref{thm:6convergenceforveechOmega} covers their normalization, as well as every norm on $\R^4$. The exact constants \eqref{eq:3convergenceforveech} can be derived from Burrin-Fairchild \cite{BurFai}.
\end{rmk}

While $\kappa_\Vc$ (resp. $\phi_\Vc$) in Theorem \ref{thm:3convergenceforveech} can be computed as a weighted sum of integrals on $\R^2$ depending on the length data of the associated Veech curve (resp. an enumeration of $\SL_2(\Vc)$-orbits in $\R^4$), it is in general difficult to say anything more concrete. In the special case of Veech surfaces arising from unfolding regular polygons, these terms can be computed more explicitly using work of Fairchild \cite{Fai}; this is done in Section \ref{subsec:heckegroups}. However, the constants are always arithmetic, in the following sense.

\begin{cor}
    For every Veech curve $\Vc$ and every $A \ge 0$, one has
    \begin{equation*}
        c_A(\Vc) \in \kay_\Vc(\pi,\Vol(\Vc^1)),
    \end{equation*}
    where $\kay_\Vc$ is the trace field of $\Vc$.
\end{cor}

\begin{rmk}
    It is important to note here that we are viewing the collection of holonomy vectors $\Hol_\omega$ as a \textit{set} as opposed to a \textit{multiset}. This simplifies the combinatorics in our proofs considerably, although the corresponding constants $c_A$ for the \textit{multiset} $\Hol_\omega$ are easily recoverable from the methods of this paper; see Section \ref{subsec:billiards} for an example.
\end{rmk}

\subsection{A brief history} Siegel-Veech constants associated to strata and their subvarieties have been extensively studied beginning with work of Veech \cite{Vee} (the case of $\SL_2(\R)/\SL_2(\Z)$ was considered first by Siegel \cite{Sie}). Computations in the one-dimensional case have been carried out by a number of authors, using a variety of techniques. We refer the interested reader to \cite{Sch}, \cite{EskMasZor}, \cite{EskMasSch}, \cite{Sch1}, \cite{Gou}, \cite{CheMolSauZag}, \cite{DelGouZogZor}, \cite{AnaMulNor}, \cite{AulCalMatSalSch}, \cite{Ngu} for a (long but) incomplete set of references. Many of these results involve techniques from algebraic geometry, and in particular interpret Siegel-Veech constants via boundary degenerations and intersection theory on the Deligne-Mumford compactification $\overline{\Ms}_{g,n}$. On the other hand, pair Siegel-Veech constants have only been considered more recently in \cite{AthCheMas}, \cite{AthFaiMas}; for Veech curves these counts have also been studied in \cite{Fai}, \cite{BurFai}. It is at the moment unclear whether or not the pair constants $c_A$ admit such interpretations, since for $A > 0$ the ``thick part'' of the moduli space will contribute nontrivially, and in fact degeneration to the boundary seems to force $A \downarrow 0$.

\subsection{Notation} We write $X \lesssim Y$ or $X = O(Y)$ if $X \le CY$ for some constant $C > 0$, and $X \asymp Y$ means $X \lesssim Y \lesssim X$. The notation $X \sim Y$ means $X/Y \to 1$, and $X= o(Y)$ means $X/Y \to 0$. We write $X \ll Y$ if $X \le cY$ for some $c > 0$ which is assumed to be sufficiently small.  The Euclidean ball of radius $r$ in $\R^m$ is denoted by $\B_r^m$.

\subsection{Acknowledgments} Research partially funded by NSF RTG grant No. DMS-2230900. This project grew out of an ongoing collaboration with Samantha Fairchild, Kanishka Katipearachchi, and Howard Masur; I would like to thank them for helpful discussions and interest in this project. I would also like to thank my advisor Paul Apisa for discussions on monodromy, leading to a simplification in the proof of Proposition \ref{prop:combinatorialclassification}.

\subsection{AI Disclosure} All figures were produced with the help of Codex. The proof of the auxiliary Lemma \ref{lemma:aasgoodpartitions} was found by ChatGPT, and was subsequently simplified by the author with AI assistance. ChatGPT was also used for proofreading, and in doing so located the following minor (mathematical) errors: a few normalization issues in Section \ref{sec:veech}, and a uniformity error in the proof of Lemma \ref{lemma:lightconeconvergence}. Besides these use cases, all ideas, arguments, examples, and writing are solely due to the author.

\section{Marked tori}

In this section, our primary goal is to establish Theorem \ref{thm:1convergenceforbranchedcovers} for the case of marked tori; this will form the basis of the proof for general $\Ec$.

\subsection{Holonomy}

Throughout, we identify $\C \cong \R^2$. Let $\Hc(0,0^{n-1})$ denote the space of tuples $(E,\Sigma)$, where $E = \C/\Lambda$ is a complex elliptic curve and $\Sigma = \{\xi_1,...,\xi_n\}$ is a set of $n$ distinct points, with $\xi_1 = 0$ being the origin in $E$. This space can be viewed in two ways: the first is as the \textit{homogeneous space}
\begin{equation*}
    \Hc(0,0^{n-1}) = \bigl(\GL_2^+(\R) \ltimes \R^{2(n-1)}\bigr) \big/\bigl(\SL_2(\Z) \ltimes \Z^{2(n-1)}\bigr),
\end{equation*}
and the second is as a space of (marked) translation surfaces. We will make use of both of these perspectives throughout.

A \textit{saddle connection} $\gamma$ on a translation surface $\omega$ is a flat geodesic segment joining cone or marked points of $\omega$, without any distinguished points in its interior; its \textit{holonomy} is $\hol_\omega(\gamma) = \int_\gamma \omega$. The set $\Hol_\omega$ of holonomy vectors of $\omega$ is then
\begin{equation*}
    \Hol_\omega = \left\{ \hol_\omega(\gamma) : \text{$\gamma$ is a saddle connection on $\omega$} \right\}.
\end{equation*}

We say $c_1,...,c_k \in \R$ are \textit{quadratically independent} if $P(c_1,...,c_k) \ne 0$ for any nonzero polynomial $P \in \Q[x_1,...,x_k]$ of degree at most 2.

\begin{lemma}\label{lemma:quadraticallyindependentholonomy}
    Consider a marked torus $\omega = (\T^2,\Sigma) \in \Hc(0,0^{n-1})$. Assume the set of coordinates $\{(\xi_i)_j : 2 \le i \le n,\, j = 1,2\}$ is quadratically independent, and let $\zeta_{ij} = \xi_j - \xi_i$. Then $\Hol_\omega = \Lambda_\Sigma$, where
    \begin{equation}\label{eq:quadraticallyindependentholonomy}
        \Lambda_{\Sigma} = \Z_{\prim}^2 \cup \bigcup_{i \ne j} \Z^2 + \zeta_{ij}.
    \end{equation}
\end{lemma}
\begin{proof}
    Suppose $\gamma$ is a saddle connection on $\omega$. We consider two cases.

    \textbf{Case I}: $\gamma$ is absolute, joining the marked point $\xi_i$ to itself. Then $\hol_\omega(\gamma) \in \Z^2$, and moreover if $\gamma$ defines a primitive curve $\hol_\omega(\gamma) \in \Z_{\prim}^2$.

    On the other hand, consider for each $(p_1,p_2) \in \Z_{\prim}^2$ the flat geodesic $\gamma(s) = (sp_1,sp_2)$. If the projection of $\gamma$ to $\T^2$ were to hit a singular point in its interior, then for some $s \in (0,1)$ and some $i$ we have
    \begin{equation*}
        s(p_1,p_2) \pm \xi_i \in \Z^2.
    \end{equation*}
    Then for some $(q_1,q_2) \in \Z^2$, we have
    \begin{equation*}
        \begin{cases}
            sp_1 \pm \xi_{i,1}+q_1 = 0 \\
            sp_2 \pm \xi_{i,2}+q_2 = 0
        \end{cases}
    \end{equation*}
    If $i=1$, then it is straightforward to see that this condition would imply $(p_1,p_2)$ is not primitive. On the other hand, if $i \ne 1$ and assuming without loss of generality that $p_1 \ne 0$, we see that
    \begin{equation*}
        \frac{p_2}{p_1}(\mp \xi_{i,1} - q_1) \pm \xi_{i,2} + q_2 = 0.
    \end{equation*}
    This is a nontrivial linear relation, which is a contradiction. Hence $\Hol_\omega^{\mathrm{abs}} = \Z_{\prim}^2$.

    \textbf{Case II}: $\gamma$ is relative, joining distinct marked points $\xi_i,\xi_j$. Then $\hol_\omega(\gamma) \in \Z^2 + \xi_j - \xi_i$; we show the reverse inclusion. Assume $v = \xi_j-\xi_i + (p_1,p_2)$, and consider the flat geodesic $\gamma(s) = \xi_i + sv$. If $\gamma$ were to hit a cone point in its interior, then there would exist some $s \in (0,1)$, some $k$ and some $(q_1,q_2) \in \Z^2$ for which
    \begin{equation*}
        \begin{cases}
            (1-s)\xi_{i,1} + s\xi_{j,1} \pm \xi_{k,1} + sp_1+q_1 = 0 \\
            (1-s)\xi_{i,2} + s\xi_{j,2} \pm \xi_{k,2} + sp_2+q_2 = 0
        \end{cases}
    \end{equation*}
    Expanding out this expression yields a nontrivial quadratic $\Q$-algebraic relation amongst the coefficients.
\end{proof}

\begin{rmk}
    See \cite{San} for a similar analysis in the case $\Hc(0,0)$. The condition of quadratic independence is much stronger than needed to establish $\Hol_\omega = \Lambda_\Sigma$; for example, the surface $\omega = (\T^2,\xi_1,\xi_2)$ with $\xi_2 = (1/2,\sqrt{2})$ does not satisfy the conditions of Lemma \ref{lemma:quadraticallyindependentholonomy}, but nonetheless satisfies \eqref{eq:quadraticallyindependentholonomy}. See Example \ref{ex:extrasubspaces}.
\end{rmk}

\subsection{Quadratic forms on $\R^4$}

Since we would like to view holonomy vectors in $\C$ as column vectors in $\R^2$, we consider the slightly nonstandard identification $\C^2 \cong \Mat_{2 \times 2}(\R)$, with coordinates
\begin{equation*}
    \Mat_{2 \times 2}(\R) = \left\{ \begin{bmatrix}
        x_1 & x_3 \\
        x_2 & x_4
    \end{bmatrix} : (x_1,x_2), (x_3,x_4) \in \C \right\}.
\end{equation*}
We let $Q(x) = x_1x_4-x_2x_3$ denote the standard determinant form on $\Mat_{2 \times 2}(\R)$; this is a signature $(2,2)$-form of discriminant 1. For $\eta \in \R^4$, define the associated inhomogeneous form $Q_\eta$ by $Q_\eta(x) = Q(x+\eta)$.

Let $\Starbodies$ denote the class of star-bodies in $\R^4$, parameterized by $\Omega = \{x \in \R^4 : \norm{x} \le \nu(x/\norm{x}) \}$ for some continuous $\nu \colon \S^3 \to (0,\infty)$. Given $\Omega \in \Starbodies$, we define
\begin{align*}
    N_{\Omega,A}(\omega,R) &= \left| \left\{ (z,w) \in \Hol_\omega \times \Hol_\omega : (z,w) \in R\Omega,\ |z \wedge w| \le A \right\} \right| \\
    c_{\Omega,A}(\omega) &= \lim_{R \to \infty} R^{-2}N_{\Omega,A}(\omega,R).
\end{align*}
Note that $N_{\B_1^4,A} = N_A$ and $c_{\B_1^4,A} = c_A$. We also define the integrals
\begin{equation}
    \kappa_\Omega^\pm = \int_{\R^2} \ind_\Omega(x,\pm x)\dx, \quad \quad \kappa_\Omega = \kappa_\Omega^+ + \kappa_\Omega^-, \quad \quad \lambda_\Omega = \int_{\Omega \cap L} \frac{d\Hs^3}{\norm{\nabla Q}}, \label{eq:def-kappalambdaOmega}
\end{equation}
where $L = \{Q = 0\}$ is the light cone, and $\Hs^s$ is the standard $s$-dimensional Hausdorff measure on $\R^4$. By \cite[Lemma 3.8]{EskMarMoz}, one has
\begin{equation*}
    \Vol(\{ x \in \R^4 : x \in R\Omega,\, a < Q_\eta(x) < b \}) \sim \lambda_\Omega(b-a)R^2.
\end{equation*}

For $t \ne 0$, let $\Sigma_t = \{Q = t\}$. Note that $\Sigma_t \to L$ in the Hausdorff topology as $t \to 0$; the following lemma shows that the corresponding surface measures converge along well-behaved sets.

\begin{lemma}\label{lemma:lightconeconvergence}
    Suppose $\Omega \in \Starbodies$. Then
    \begin{equation*}
        \lim_{t \to 0} \int_{\Omega \cap \Sigma_t} \frac{\dH^3}{\norm{\nabla Q}} = \lambda_\Omega.
    \end{equation*}
\end{lemma}
\begin{proof}
    Set $\ebf_\theta = (\cos \theta,\sin \theta)$. Make the change of variables
    \begin{equation*}
        (z,w) = \bigl( r\ebf_\theta, s\ebf_\theta + tr^{-1}\ebf_{\theta+\pi/2} \bigr).
    \end{equation*}
    In these coordinates $\frac{\dH^3}{\norm{\nabla Q}} = \dr\ds\dtheta$, so
    \begin{align*}
        \int_{\Omega \cap \Sigma_t} \frac{\dH^3}{\norm{\nabla Q}} &= \int_{0}^\infty \int_{-\infty}^\infty \int_0^{2\pi} \ind_\Omega \bigl( r\ebf_\theta, s\ebf_\theta + tr^{-1}\ebf_{\theta+\pi/2} \bigr)\dtheta \ds\dr \\
        &= \int_{0}^{\diam \Omega} \int_{-\diam \Omega}^{\diam \Omega} \int_0^{2\pi} \ind_\Omega \bigl( r\ebf_\theta, s\ebf_\theta + tr^{-1}\ebf_{\theta+\pi/2} \bigr)\dtheta \ds\dr.
    \end{align*}
    Since $\del \Omega \cap L$ is a null set, we may apply dominated convergence to obtain the desired result.
\end{proof}

\begin{ex}\label{ex:kappalambdacomputation1}
    A computation shows
    \begin{equation*}
        \kappa_{\B_1^4} = \pi, \quad \quad \lambda_{\B_1^4} = \pi^2, \quad \quad \kappa_{\B_1^2 \times \B_1^2} = 2\pi, \quad \quad \lambda_{\B_1^2 \times \B_1^2} = 4\pi.
    \end{equation*}
    For the domain $\Omega_{\AFM} = \{(z,w) \in \C^2 : |w| \le |z| \le 1\}$, we have $\kappa_{\Omega_{\AFM}} = \lambda_{\Omega_{\AFM}} = 2\pi$. On the other hand, strictly speaking $\Omega_{\AFM} \notin \Starbodies$, but our results may be applied by approximating by (unions and differences of) elements of $\Starbodies$ and noting that $\int_{\del \Omega_{\AFM} \cap L} \frac{\dH^3}{\norm{\nabla Q}} = 0$.
\end{ex}

Given a vector $\xi \in \R^m$, we say $\xi$ is \textit{Diophantine} if there exist $C = C(\xi) > 0$ and $\mu = \mu(\xi) > 0$ so that for all $Q \ge 1$ and all $\qbf = (p_1/q_1,...,p_m/q_m) \in \Q^m$ with $1 \le q_i \le Q$, we have
\begin{equation*}
    \norm{\xi - \qbf} \ge \frac{C}{Q^\mu}.
\end{equation*}
It is known that the set of $\xi \in \R^m$ which are not Diophantine has Hausdorff dimension zero.

If $V \subseteq \R^4$ is a two-dimensional rational nullspace of $Q$ and $\eta \in V + x_\eta$ for some $x_\eta \in \Z^4$, we say $V - x_\eta$ is an \textit{exceptional subspace} for $Q_\eta$. For $\eta \in \R^4$, $\Omega \in \Starbodies$, and $(a,b) \subseteq \R$, we define the counting function
\begin{equation*}
    \Tilde{N}_{\eta,\Omega}(a,b,R) = \left| \left\{ x \in \Z^4 : \substack{\text{$x$ is not in a $Q_\eta$-exceptional subspace} \\ x \in R\Omega,\ Q_\eta(x) \in (a,b)} \right\} \right|.
\end{equation*}
The following is a special case of \cite{MarMoh}.

\begin{thm}[\cite{MarMoh}, Theorem 1.9]\label{thm:MarMoh}
    Suppose $\eta \in \R^4$ is Diophantine. Then for any $\Omega \in \Starbodies$ and any interval $(a,b) \subseteq \R$, we have
    \begin{equation*}
        \Tilde{N}_{\eta,\Omega}(a,b,R) \sim \lambda_\Omega(b-a)R^2.
    \end{equation*}
\end{thm}

The connection between Theorems \ref{thm:MarMoh} and \ref{thm:1convergenceforbranchedcovers} lies in the observation that, if $\omega = (\T^2,\Sigma) \in \Hc(0,0^{n-1})$ and $z,w \in \Hol_\omega$, then we may express
\begin{equation*}
    z = \begin{bmatrix}
        p_1 \\
        q_1
    \end{bmatrix} + \zeta_{ij}, \quad \quad w = \begin{bmatrix}
        p_2 \\
        q_2
    \end{bmatrix} + \zeta_{kl}
\end{equation*}
for some $1 \le i,j,k,l \le n$ and some $p_1,p_2,q_1,q_2 \in \Z$; then
\begin{equation*}
    z \wedge w = \det \begin{bmatrix}
        p_1 + (\zeta_{ij})_1 & p_2 + (\zeta_{kl})_1 \\
        q_1 + (\zeta_{ij})_2 & q_2 + (\zeta_{kl})_2 \\
    \end{bmatrix}
    = Q_\eta(p_1,q_1,p_2,q_2),
\end{equation*}
where $\eta = (\zeta_{ij},\zeta_{kl})$. Thus, pair counting on $\Hc(0,0^{n-1})$ can be recast as understanding the distribution of inhomogeneous quadratic forms along (subsets of) integer points.

\subsection{The Diophantine case}\label{subsec:diophantine}

In this section, we establish a generalization of Theorem \ref{thm:1convergenceforbranchedcovers} when $\Ec = \Hc(0,0^{n-1})$.

\begin{defn}[Diophantine marked tori]\label{defn:diophantinemarkedtori}
    We say $\omega = (\T^2,\Sigma)$ is \textit{weakly Diophantine} if the differences $\zeta_{ij} = \xi_j-\xi_i$ are distinct Diophantine vectors in $\R^2$. It is \textit{Diophantine} if it is weakly Diophantine, and additionally satisfies the quadratic independence condition of Lemma \ref{lemma:quadraticallyindependentholonomy}.
\end{defn}

Let us define
\begin{equation}\label{eq:def-FOmegadnA}
    F_{\Omega,d,n}(A) = \kappa_\Omega \biggl[ \frac{d}{\zeta(2)} + dn(n-1) \biggr] + \lambda_\Omega\biggl[ \frac{2d}{\zeta(2)}\Phi_3(dA) + d^2\biggl(\frac{4n(n-1)}{\zeta(2)} + 2n^2(n-1)^2 \biggr)A \biggr].
\end{equation}
The following theorem, combined with Example \ref{ex:kappalambdacomputation1}, generalizes Theorem \ref{thm:1convergenceforbranchedcovers}.

\begin{thm}\label{thm:4convergencefordiophantinetori}
    Suppose $\omega = (\T^2,\Sigma)$ is Diophantine. Then for every $\Omega \in \Starbodies$ and $A \ge 0$,
    \begin{equation*}
        c_{\Omega,A}(\omega) = F_{\Omega,1,n}(A).
    \end{equation*}
\end{thm}

\begin{rmk}
    Note that $F_{\Omega,d,n}(A)$ is essentially a polynomial in $d,n,A$ (up to the $\Phi_3$-correction factor). Moreover, the coefficients of $F_{\Omega,d,n}$ lie in the field $\Q(\Vol(\Hc^1(0)),\kappa_\Omega,\lambda_\Omega)$, and asymptotically
\begin{equation}\label{eq:largeAlimit}
    \lim_{A \to \infty} \frac{F_{\Omega,d,n}(A)}{A} = 2\lambda_\Omega d^2 \left( n(n-1) + \frac{1}{\Vol(\Hc^1(0))} \right)^2.
\end{equation}
\end{rmk}

\begin{rmk}
    One can define $c_{\Omega,I}(\omega)$ for an arbitrary bounded interval $I \subseteq \R$, as opposed to restricting to $I = [-A,A]$, but due to the symmetry of $\Hol_\omega$ under multiplication by $-1$ this yields no additional information.
\end{rmk}

To begin the proof of Theorem \ref{thm:4convergencefordiophantinetori}, we decompose the counting function $N_{\Omega,A}(\omega,R)$ into \textit{absolute}, \textit{relative}, and \textit{mixed} pieces:
\begin{equation}\label{eq:NAdecomp}
    N_{\Omega,A}(\omega,R) = N_{\Omega,A}^{\abs}(\omega,R) + \sum_{\substack{1 \le i,j,k,l \le n \\i\ne j, k \ne l}}N_{\Omega,A}^{\rel,ijkl}(\omega,R) + \sum_{\substack{1 \le i,j \le n \\ i \ne j}} N_{\Omega,A}^{\mix,ij}(\omega,R),
\end{equation}
where the subterms in \eqref{eq:NAdecomp} are defined by
\begin{align*}
    N_{\Omega,A}^{\abs}(\omega,R) &= \left| \left\{ (z,w) \in \Z_{\prim}^2 \times \Z_{\prim}^2 : (z,w) \in R\Omega,|z \wedge w| \le A \right\} \right|, \\
    N_{\Omega,A}^{\rel,ijkl}(\omega,R) &= \left| \left\{ (z,w) \in (\Z^2 + \zeta_{ij}) \times (\Z^2 + \zeta_{kl}) : (z,w) \in R\Omega,|z \wedge w| \le A \right\} \right|, \\
    N_{\Omega,A}^{\mix,ij}(\omega,R) &= \left| \left\{ (z,w) \in (\Z^2 + \zeta_{ij}) \times \Z_{\prim}^2 \cup \Z_{\prim}^2 \times (\Z^2 + \zeta_{ij}) : (z,w) \in R\Omega,|z \wedge w| \le A \right\} \right|.
\end{align*}
Define $c_{\Omega,A}^*(\omega)$ for each of the above counting functions analogously to \eqref{eq:def-cA}.

\subsubsection{Absolute pairs}

\begin{lemma}\label{lemma:diophantineabscomputation}
    Under the assumptions of Theorem \ref{thm:4convergencefordiophantinetori}, we have
    \begin{equation*}
        c_{\Omega,A}^{\abs}(\omega) = \frac{\kappa_\Omega}{\zeta(2)} + \frac{2\lambda_\Omega}{\zeta(2)}\Phi_3(A).
    \end{equation*}
\end{lemma}
\begin{proof}
    Note that if $z,w \in \Z_{\prim}^2$, then $z \wedge w \in \Z$, and so our determinant can take only discrete values. Therefore $N_{\Omega,A}^{\abs}$ further decomposes as a finite sum over determinants in $[0,A] \cap \Z$.

    We begin with the case $z \wedge w = 0$. Then $z,w$ are parallel, and standard counting estimates imply
    \begin{equation}\label{eq:absolute0parallelcontribution}
        |\{z,w \in \Z_{\prim}^2 : (z,w) \in R\Omega, z \wedge w =0 \}| \sim \frac{\kappa_\Omega}{\zeta(2)}R^2.
    \end{equation}
    
    Now, consider the case of nonzero determinant. By Schmidt \cite[Lemma 6]{Schm}, for each $n \in \Z$ nonzero, there are precisely $\phi(n)$ pairs $(z,w) \in \Z_{\prim}^2 \times \Z_{\prim}^2$ with $z \wedge w = n$, up to the $\SL_2(\Z)$-action. Each of these classes has a representative of size $O_n(1)$, so for $n$ fixed and $R \gg_n 1$,
    \begin{equation*}
        |\{z,w \in \Z_{\prim}^2 : (z,w) \in R\Omega, |z \wedge w| =n \}| \sim 2\phi(n)\sum_{\substack{\gamma \in \SL_2(\Z) \\ \gamma(z_0,w_0) \in R\Omega}} 1,
    \end{equation*}
    where $(z_0,w_0)$ is some fixed pair of primitive holonomy vectors with determinant $n$. Our goal is to convert the sum on the right-hand side of this equation into a count of lattice points on a certain principal homogeneous space. To do this, recall $\Sigma_t = \{x \in \R^4 : Q(x) = t\}$, and let $\Omega_{t,R} = R\Omega \cap \Sigma_t$. For $t \ne 0$, $\Sigma_t$ is a principal homogeneous space for $\SL_2(\R)$, and hence for a fixed choice $\dg$ of Haar measure on $\SL_2(\R)$ we obtain by pushforward a unique measure $\sigma_t$ on $\Sigma_t$. Under the normalization $\int_{\SL_2(\R)/\SL_2(\Z)} \dg = \zeta(2)$ we have $\sigma_t = \frac{\dH^3}{t\norm{\nabla Q}}$. By the assumptions on the domain $\Omega$, we see that $(\Omega_{n,R})_{R \gg_n 1}$ forms a \textit{well-rounded family} in the sense of Eskin-McMullen \cite{EskMcM}, and so by \cite[Theorem 1.4]{EskMcM} (see also Gorodnik-Nevo \cite{GorNev}),
    \begin{equation*}
        \sum_{\substack{\gamma \in \SL_2(\Z) \\ \gamma(z_0,w_0) \in R\Omega}} 1 \sim \frac{\sigma_n(\Omega_{n,R})}{\Covol(\SL_2(\Z))} = \frac{\sigma_n(\Omega_{n,R})}{\zeta(2)}.
    \end{equation*}
    To find an asymptotic for $\sigma_n(\Omega_{n,R})$, note
    \begin{equation*}
        \sigma_n(\Omega_{n,R}) = n^{-1}\int_{R\Omega \cap \Sigma_n} \frac{\dH^3}{\norm{\nabla Q}} = n^{-1}R^2 \int_{\Omega \cap \Sigma_{n/R^2}} \frac{\dH^3}{\norm{\nabla Q}}.
    \end{equation*}
    Combining this with Lemma \ref{lemma:lightconeconvergence}, we obtain for $n \in \N$ that
    \begin{equation}\label{eq:absolutencontribution}
        |\{z,w \in \Z_{\prim}^2 : (z,w) \in R\Omega, |z \wedge w| =n \}| \sim \frac{2\lambda_\Omega \phi(n)}{n\zeta(2)}R^2.
    \end{equation}
    Summing over each $n \in [0,A] \cap \Z$ yields the desired result.
\end{proof}

\begin{rmk}
    For $\Omega = \B_1^4$, Lemma \ref{lemma:diophantineabscomputation} was proven implicitly by Schmidt \cite{Schm} (see also \cite[Section 5]{Fai}).
\end{rmk}

\subsubsection{Relative pairs}

\begin{lemma}\label{lemma:diophantinerelcomputation}
    Under the assumptions of Theorem \ref{thm:4convergencefordiophantinetori}, we have
    \begin{equation*}
        c_{\Omega,A}^{\rel,ijkl}(\omega) = \left\{\begin{array}{cc}
            2\lambda_\Omega A + \kappa_\Omega^+, & (i,j) = (k,l), \\
            2\lambda_\Omega A + \kappa_\Omega^-, & (i,j) = (l,k), \\
            2\lambda_\Omega A, & \text{otherwise}
        \end{array}\right.
    \end{equation*}
\end{lemma}
\begin{proof}
    Note that for fixed $\omega$, we have as $R \to \infty$ the asymptotic
    \begin{equation*}
        N_{\Omega,A}^{\rel,ijkl}(\omega,R) \sim \left| \left\{ x \in \Z^4 : x \in R\Omega,\ Q_{\eta_{ijkl}}(x) \in [-A,A] \right\} \right|,
    \end{equation*}
    where $\eta_{ijkl} = (\zeta_{ij},\zeta_{kl})$. The quadratic independence condition on $\omega$ implies $Q_{\eta_{ijkl}}$ has no exceptional subspaces unless $(i,j) = (k,l)$ or $(i,j)=(l,k)$. If neither of these cases hold,
    \begin{equation*}
        \left| \left\{ x \in \Z^4 : x \in R\Omega,\ Q_{\eta_{ijkl}}(x) \in [-A,A] \right\} \right| \sim \Tilde{N}_{\eta_{ijkl},\Omega}(-A,A,R) \sim 2\lambda_\Omega AR^2,
    \end{equation*}
    where the first asymptotic follows from $\{Q_{\eta_{ijkl}} = A\}$ not being an exceptional subspace, and the second asymptotic follows from Theorem \ref{thm:MarMoh}; the Diophantine assumption on $\omega$ is precisely the condition that $\eta_{ijkl}$ is a Diophantine vector for every $i,j,k,l$.

    When $(i,j) = (k,l)$ (resp. $(i,j) = (l,k)$), there is an exceptional subspace $\{(z,z) : z \in \C \}$ (resp. $\{(z,-z) : z \in \C \}$). On these 2-planes,
    \begin{equation*}
        \left| \{ x \in \Z^4 : x \in R\Omega \cap \{(z,\pm z)\}\} \right| \sim \kappa_\Omega^\pm R^2.
    \end{equation*}
    Combining these estimates yields the desired result.
\end{proof}

\subsubsection{Mixed pairs}

Now, we turn to the case of mixed pairs. This is the trickiest case, since we must restrict to primitive integer pairs in one, but not both, of the coordinates; dealing with this requires a slightly delicate M\"obius inversion argument.

\begin{prop}\label{prop:MarMohmix}
    Suppose $\eta = (0,0,\zeta_1,\zeta_2)$ is Diophantine, and let
    \begin{equation*}
        \Tilde{N}_{\eta,\Omega}^{\mix}(a,b,R) = \left| \left\{ x \in \Z^4 : \gcd(x_1,x_2) = 1,\ x \in R\Omega,\ Q_\eta(x) \in (a,b),\ Q_\eta(x) \ne 0 \right\} \right|.
    \end{equation*}
    Then
    \begin{equation*}
        \Tilde{N}_{\eta,\Omega}^{\mix}(a,b,R) \sim \zeta(2)^{-1}\lambda_\Omega(b-a)R^2.
    \end{equation*}
\end{prop}

We note that Proposition \ref{prop:MarMohmix} would follow immediately from \cite{MarMoh} by inclusion-exclusion assuming all involved asymptotics were \textit{uniform}, but of course this is unrealistic to expect (and likely to be false). To circumvent this, we instead perform a partial M\"obius inversion argument with a large cutoff, and show the resulting error term is of negligible size. For this, we require the homogeneity relation $Q_\eta(\ell u,v) = \ell Q_\eta(u,v)$, which will not hold for general inhomogeneous forms.

\begin{proof}[Proof of Proposition \ref{prop:MarMohmix}]
    Given $d \in \N$, consider the modified counting function
    \begin{align*}
        N_d(R) &= \left| \left\{ x \in \Z^4 : d \mid \gcd(x_1,x_2),\, x \in R\Omega,\, Q_\eta(x) \in (a,b),\, Q_\eta(x) \ne 0 \right\} \right| \\
        &= \left| \left\{ x' = (y_1,y_2,x_3,x_4) \in \Z^4 : x' \in R\Omega_d,\, Q_\eta(x') \in (a/d,b/d),\, Q_\eta(x') \ne 0 \right\} \right|,
    \end{align*}
    where $\Omega_d = \{v \in \R^4 : (dv_1,dv_2,v_3,v_4) \in \Omega\}$. Applying Theorem \ref{thm:MarMoh} alongside a standard Jacobian estimate yields
    \begin{equation}\label{eq:NdRasymp}
        N_d(R) \sim d^{-2}\lambda_\Omega(b-a)R^2.
    \end{equation}
    We note the convergence in \eqref{eq:NdRasymp} will not necessarily be uniform in $d$.

    Now, fix $H = h!$ a large factorial, and consider
    \begin{equation*}
        N^H(R) = \left| \left\{ x \in \Z^4 : \gcd(x_1,x_2,H) = 1,\, x \in R\Omega,\, Q_\eta(x) \in (a,b),\, Q_\eta(x) \ne 0 \right\} \right|.
    \end{equation*}
    That is, $N^H(R)$ counts the number of integer points $x \in \Z^4$ contributing to the count of $\Tilde{N}_{\eta,\Omega}(a,b,R)$, so that $\gcd(x_1,x_2)$ does not have a prime factor of size less than $h$. By M\"obius inversion, we have
    \begin{align}
        N^H(R) &= \sum_{\substack{x \in \Z^4 \\ x \in R\Omega \\ Q_\eta(x) \in (a,b) \setminus \{0\}}} \ind_{\gcd(x_1,x_2,H) = 1} \notag \\
        &= \sum_{\substack{x \in \Z^4 \\ x \in R\Omega \\ Q_\eta(x) \in (a,b) \setminus \{0\}}} \sum_{d \mid \gcd(x_1,x_2,H)} \mu(d) = \sum_{d \mid H} \mu(d)N_d(R). \label{eq:Mobiusinversion}
    \end{align}

    Now, if we fix $H$, we have by \eqref{eq:NdRasymp} that
    \begin{equation}\label{eq:MobiusinversionfixedHlimit}
        \lim_{R \to \infty} R^{-2} N^H(R) = (1+o_{H \to \infty}(1))\zeta(2)^{-1}\lambda_\Omega(b-a).
    \end{equation}
    Moreover, we have the upper bound
    \begin{equation}\label{eq:MobiusinversionmixtoH}
        \left|\Tilde{N}_{\eta,\Omega}^{\mix}(a,b,R) - N^H(R) \right| \le \sum_{\substack{(u,v) \in \Z_{\prim}^2 \times \Z^2 \\ Q_\eta(u,v) \ne 0 \\ (u,v) \in R\B_{C_2}^4}} \sum_{\substack{\ell \ge 2 \\ \gcd(\ell,H) = 1 \\ (\ell u,v) \in R\B_{C_2}^4}} \ind_{(-C_1,C_1)}(Q_\eta(\ell u,v)),
    \end{equation}
    where $C_1 = \max\{|a|,|b|\}$ and $C_2 = \diam(\Omega)$. Setting
    \begin{equation*}
        X(u,v,R) = \left\lfloor \min \left\{ \frac{C_1}{|u \wedge (v+\zeta)|}, \frac{C_2R}{|u|} \right\} \right\rfloor,
    \end{equation*}
    then we see that the sum over $\ell$ in \eqref{eq:MobiusinversionmixtoH} is restricted to the range $2 \le \ell \le X(u,v,R)$. A classical sieve-theoretic argument (see e.g. \cite{Ten}) shows
    \begin{equation*}
        \Phi(X,H) = \left| \left\{ 2 \le \ell \le X : \gcd(\ell,H) = 1 \right\} \right| \lesssim \frac{X}{\log h}
    \end{equation*}
    uniformly in $X \ge 0$ and $h \ge 2$.\footnote{The function $\Phi$ is (up to adding 1) the count of $h$-rough numbers of size at most $X$. The stated uniform inequality is usually restricted to the range $2 \le h \le X$, but of course for $h > X$ this count is trivially zero.} Therefore, we may upper bound \eqref{eq:MobiusinversionmixtoH} by
    \begin{equation}\label{eq:MobiusinversionHtologh}
        \eqref{eq:MobiusinversionmixtoH} \le \sum_{\substack{(u,v) \in \Z_{\prim}^2 \times \Z^2 \\ Q_\eta(u,v) \ne 0 \\ (u,v) \in R\B_{C_2}^4}} \Phi(X(u,v,R),H) \lesssim \frac{1}{\log h} \sum_{\substack{(u,v) \in \Z_{\prim}^2 \times \Z^2 \\ Q_\eta(u,v) \ne 0 \\ (u,v) \in R\B_{C_2}^4}} X(u,v,R).
    \end{equation}
    On the other hand, we observe
    \begin{align*}
        \sum_{\substack{(u,v) \in \Z_{\prim}^2 \times \Z^2 \\ Q_\eta(u,v) \ne 0 \\ (u,v) \in R\B_{C_2}^4}} X(u,v,R) &\le \sum_{\substack{(u,v) \in \Z_{\prim}^2 \times \Z^2 \\ Q_\eta(u,v) \ne 0 \\ (u,v) \in R\B_{C_2}^4}} \left| \left\{\ell \ge 1 : (\ell u,v) \in R\B_{2C_2}^4,\, Q_\eta(\ell u,v) \in (-2C_1,2C_1) \right\} \right| \\
        &= \Tilde{N}_{\eta,\B_{2C_2}^4}(-2C_1,2C_1,R),
    \end{align*}
    and so we have the bound
    \begin{equation}\label{eq:uniformboundlogh}
        \frac{1}{\log h} \sum_{\substack{(u,v) \in \Z_{\prim}^2 \times \Z^2 \\ Q_\eta(u,v) \ne 0 \\ (u,v) \in R\B_{C_2}^4}} X(u,v,R) \le \frac{1}{\log h}\Tilde{N}_{\eta,\B_{2C_2}^4}(-2C_1,2C_1,R).
    \end{equation}
    Theorem \ref{thm:MarMoh} implies that the right-hand side of \eqref{eq:uniformboundlogh} is at most $\frac{1}{\log h}C_3R^2$ for some uniform constant $C_3 = O_{\eta,C_1,C_2}(1)$, independent of $R,H$. Combining \eqref{eq:MobiusinversionmixtoH}, \eqref{eq:MobiusinversionHtologh}, \eqref{eq:uniformboundlogh} yields
    \begin{equation}\label{eq:boundbylogh}
        \limsup_{R \to \infty} R^{-2}\left| \Tilde{N}_{\eta,\Omega}^{\mix}(a,b,R) - N^H(R) \right| \lesssim \frac{C_3}{\log h}.
    \end{equation}

    Finally, we bound
    \begin{align*}
        \limsup_{R \to \infty} \left| R^{-2}\Tilde{N}_{\eta,\Omega}^{\mix}(a,b,R) - \frac{\lambda_\Omega(b-a)}{\zeta(2)} \right| &\le \limsup_{h \to \infty} \limsup_{R \to \infty} R^{-2}\left| \Tilde{N}_{\eta,\Omega}^{\mix}(a,b,R) - N^H(R) \right| \\
        &\quad \quad + \limsup_{h \to \infty} \limsup_{R \to \infty} \left| R^{-2}N^H(R) - \frac{\lambda_\Omega(b-a)}{\zeta(2)} \right|.
    \end{align*}
    The first term in the above inequality is zero by \eqref{eq:boundbylogh}, and the second term is zero by \eqref{eq:MobiusinversionfixedHlimit}.
\end{proof}

As previously noted, the homogeneity relation $Q_\eta(\ell u,v) = \ell Q_\eta(u,v)$ was used nontrivially in multiple places throughout the proof, and in particular both \eqref{eq:MobiusinversionmixtoH} and \eqref{eq:MobiusinversionHtologh} required restricting to the \textit{particular} interval $2 \le \ell \le X$; clearly an inequality of the form $\Phi(X,H) \lesssim \frac{X}{\log h}$ cannot hold when $\ell$ is sampled along an arbitrary set of $\sim X$ integers.

\begin{que}\label{que:MarMohmix}
    Does the conclusion of Proposition \ref{prop:MarMohmix} hold for arbitrary inhomogeneous forms?
\end{que}

To explain the difficulty of Question \ref{que:MarMohmix}, note that the Siegel transform $\Theta^{\mix}$ defined by
\begin{equation*}
    \Theta^{\mix}f(g,\eta) = \sum_{x \in (\Z_{\prim}^2 \times \Z^2) + \eta}f(gx)
\end{equation*}
is no longer invariant under the full $\SO(Q)$-action; contrast this with $\Z_{\prim}^4$.

\begin{lemma}\label{lemma:diophantinemixcontribution}
    Under the assumptions of Theorem \ref{thm:4convergencefordiophantinetori}, we have
    \begin{equation*}
        c_{\Omega,A}^{\mix,ij}(\omega) = \frac{4\lambda_\Omega}{\zeta(2)}A.
    \end{equation*}
\end{lemma}
\begin{proof}
    Suppose first that we are in the case $\Z_{\prim}^2 \times (\Z^2 + \zeta_{ij})$. We follow the strategy of Lemma \ref{lemma:diophantinerelcomputation} with the form $Q_{\eta_{ij}'}$, where $\eta_{ij}' = (0,\zeta_{ij})$. Note that $Q_{\eta_{ij}'}$ does have an exceptional subspace $V_{12} = \{x_1=x_2=0\}$; however, since $z$ is primitive, we always have $(z,w) \notin V_{12}$, and so by Proposition \ref{prop:MarMohmix} we have the asymptotic
    \begin{equation*}
        \left| \left\{ (z,w) \in \Z_{\prim}^2 \times (\Z^2 + \zeta_{ij}) : (z,w) \in R\Omega,\, |z \wedge w| \le A \right\} \right| \sim \Tilde{N}_{\eta_{ij}',\Omega}^{\mix}(-A,A,R) \sim \frac{2\lambda_\Omega A}{\zeta(2)}R^2.
    \end{equation*}
    For the case $(\Z^2 + \zeta_{ij}) \times \Z_{\prim}^2$, replace $\Omega$ with $\Bar{\Omega} = \{(z,w) \in \C^2 : (w,z) \in \Omega\}$ and argue as above, noting that $\lambda_{\Bar{\Omega}} = \lambda_\Omega$.
\end{proof}

\subsubsection{The proof of Theorem \ref{thm:4convergencefordiophantinetori}}

\begin{proof}[Proof of Theorem \ref{thm:4convergencefordiophantinetori}]
    Using the decomposition \eqref{eq:NAdecomp}, we see that
    \begin{equation}\label{eq:cOmegaAdecomp}
        c_{\Omega,A}(\omega) = c_{\Omega,A}^{\abs}(\omega) + \sum_{\substack{1 \le i,j,k,l \le n \\ i \ne j, k \ne l}} c_{\Omega,A}^{\rel,ijkl}(\omega) + \sum_{\substack{1 \le i,j \le n \\ i \ne j}} c_{\Omega,A}^{\mix,ij}(\omega),
    \end{equation}
    assuming each of the terms is well-defined. By Lemmas \ref{lemma:diophantineabscomputation}, \ref{lemma:diophantinerelcomputation}, \ref{lemma:diophantinemixcontribution}, the limits on the right-hand side exist, and moreover
    \begin{equation*}
        \eqref{eq:cOmegaAdecomp} = \frac{\kappa_\Omega}{\zeta(2)} + \frac{2\lambda_\Omega}{\zeta(2)}\Phi_3(A) + 2n^2(n-1)^2\lambda_\Omega A + \kappa_\Omega n(n-1) + n(n-1)\frac{4\lambda_\Omega}{\zeta(2)}A.
    \end{equation*}
    Rearranging yields the desired result.
\end{proof}

Note that the set of tuples $(\T^2,\Sigma)$ is a null set (in fact, codimension 3) in $\Hc^1(0,0^{n-1})$, and so Theorem \ref{thm:4convergencefordiophantinetori} does not by itself establish Theorem \ref{thm:1convergenceforbranchedcovers} in the case $\Ec = \Hc(0,0^{n-1})$; however, an easy argument allows us to extend to the entire space of marked tori.

Suppose $\omega \in \Hc^1(0,0^{n-1})$ has underlying torus $E \in \Hc^1(0)$. Then there exists some $g \in \SL_2(\R)$ so that $g\omega$ has underlying torus $\T^2$; we say $\omega$ is \textit{Diophantine} if $g\omega$ is. (Note that the matrix $g \in \SL_2(\R)$ is not unique -- all such $g$ are of the form $g = \gamma g_0$ for $\gamma \in \SL_2(\Z)$ and some fixed $g_0 \in \SL_2(\R)$ -- but the definition is independent of this choice).

\begin{cor}\label{cor:convergencefordiophantinetori}
    Suppose $\omega \in \Hc^1(0,0^{n-1})$ is Diophantine. Then for all $\Omega \in \Starbodies$ and all $A \ge 0$,
    \begin{equation*}
        c_{\Omega,A}(\omega) = F_{\Omega,1,n}(A).
    \end{equation*}
    In particular, the set of surfaces in $\Hc^1(0,0^{n-1})$ for which this equality does not hold is (Hausdorff) codimension at least 1.
\end{cor}
\begin{proof}
    Note that $\Hc^1(0,0^{n-1})$ is an $\SL_2(\R)$-equivariant fiber bundle over $\Hc^1(0)$, where the fiber over $\C/\Lambda \in \Hc^1(0)$ is identified with $(\C/\Lambda)^{n-1}$:
    \begin{equation*}
        \begin{tikzcd}
            (\C/\Lambda)^{n-1} \ar[r, hookrightarrow] & \Hc^1(0,0^{n-1}) \ar[d] \\
            & \Hc^1(0)
        \end{tikzcd}
    \end{equation*}
    $\SL_2(\R)$ acts transitively on $\Hc^1(0)$, and so for each $\C/\Lambda$ there exists $g \in \SL_2(\R)$ with $g\C/\Lambda = \T^2$. Moreover, the diagonal action of $\SL_2(\R)$ on $\C \times \C$ preserves the integrals $\kappa_\Omega,\lambda_\Omega$. To obtain the codimension result, simply note that the transition maps between fibers are locally Lipschitz, and use the standard inequality $\dim_H(F \times G) \le \dim_H F + \dim_B G$.
\end{proof}

\begin{cor}\label{cor:cAH0n}
    For all $\Omega \in \Starbodies$ and all $A \ge 0$, we have $c_{\Omega,A}(\Hc(0,0^{n-1})) = F_{\Omega,1,n}(A)$.
\end{cor}

\begin{rmk}
    In the proof of Theorem \ref{thm:4convergencefordiophantinetori}, we decomposed the counting functions $N_A(\omega,R)$ into sums over pairs $\Cc,\Dc$ of \textit{configurations} of saddle connections (see Eskin-Masur-Zorich \cite{EskMasZor} for definitions), and then computed the \textit{$(\Cc,\Dc)$-Siegel-Veech constants}
    \begin{equation*}
        c_A^{\Cc,\Dc}(\omega) = \lim_{R \to \infty} \frac{1}{R^2}\left| \left\{ (z,w) \in V_{\Cc}(\omega) \times V_\Dc(\omega) : (z,w) \in R\Omega,\, |z \wedge w| \le A \right\} \right|
    \end{equation*}
    associated to each such pair.
\end{rmk}

\begin{rmk}
    Theorem \ref{thm:4convergencefordiophantinetori} can be proven more easily for \textit{generic} surfaces $\omega \in \Hc^1(0,0^{n-1})$. A sketch of the proof is as follows: let $\hat{f}$ be the \textit{pair Siegel-Veech transform} on $\R^4$; by an argument analogous to \cite[Section 4.5.2]{AthCheMas},\footnote{\cite{AthCheMas} does not treat the case of mixed terms, but this can be easily remedied via a folding-unfolding argument.} one can show $\hat{f} \in L^1(\Hc^1(0,0^{n-1}))$, and moreover identify the associated Siegel-Veech measure $\lambda$ on $\R^4$ for which $\int_{\Hc^1(0,0^{n-1})} \hat{f}(\omega) \dmu(\omega) = \int_{\R^4} f\dlambda$. This measure decomposes into two-dimensional pieces (corresponding to parallel pairs), three-dimensional pieces (supported on the integer level sets of $Q$), and some multiple of the standard Lebesgue measure. On the $\R^4$-piece, one can apply Nevo-Stein \cite{NevSte} alongside the arguments of \cite{AthFaiMas}, whereas the lower-dimensional slices reduce to the (known) case of $\Hc(0)$.
\end{rmk}

\subsection{Holonomy vectors in general, and some examples}

Note that we always have the inclusion $\Hol_\omega \subseteq \Lambda_\Sigma$, and so Theorem \ref{thm:MarMoh} is easily seen to imply the following.

\begin{cor}\label{cor:growthrateofpairs}
    Suppose $\omega \in \Hc(0,0^{n-1})$ is weakly Diophantine. Then for each $\Omega \in \Starbodies$ and $A \ge 0$,
    \begin{equation*}
        \limsup_{R \to \infty} \frac{1}{R^2}N_{\Omega,A}(\omega;R) < \infty.
    \end{equation*}
    In particular, the set of surfaces for which this fails is codimension at least 2.
\end{cor}

We say a surface $\omega \in \Hc$ is \textit{exceptional} if $\limsup_{R \to \infty} \frac{1}{R^2}N_{\Omega,A}(\omega;R) = \infty$ for some $\Omega \in \Starbodies$ and some $A > 0$. In \cite{FaiFenKatMas}, we will show the following.

\begin{thm}[\cite{FaiFenKatMas}]\label{thm:counterexamples}
    There exists a dense $G_\delta$ codimension-2 set of surfaces $\omega \in \Hc(0,0^{n-1})$ for which
    \begin{equation*}
        N_{\Omega,A}(\omega,R) \ge \frac{R^2\log R}{\log\log R}
    \end{equation*}
    holds for all $\Omega \in \Starbodies$, all $A > 0$, and arbitrarily large $R \ge 1$.
\end{thm}

\begin{proof}[Proof of Theorem \ref{thm:2dimensiongapforcAEc} for $\Ec = \Hc(0,0^{n-1})$]
    Combine Corollary \ref{cor:growthrateofpairs} with Theorem \ref{thm:counterexamples}.
\end{proof}

We also make a remark on the necessity of the assumptions of Theorem \ref{thm:1convergenceforbranchedcovers}. The Diophantine hypotheses are clearly necessary in view of Theorem \ref{thm:counterexamples}. On the other hand, Lemma \ref{lemma:diophantineabscomputation} holds for \textit{every} surface $\omega$ satisfying the condition \eqref{eq:quadraticallyindependentholonomy}, and so the Diophantine assumptions on the marked points can be done away with in the absolute case; this observation forms the basis of the proof of Theorem \ref{thm:3convergenceforveech}. On the other hand, without any independence condition on the holonomy vectors of $\omega$, Lemma \ref{lemma:diophantineabscomputation} can fail spectacularly, as in the following example.

\begin{ex}[A marked surface with no absolute saddle connections]
    Let $E \in \Hc(0)$ have period lattice $\Lambda$, and let $E[N] = \{v \in E : Nv = 0\} = N^{-1}\Lambda/\Lambda$ denote its $N$-torsion points. Choose $\xi_1,...,\xi_n \in E$ distinct, with $\xi_1 = 0$. We define $\omega = (E,\Sigma_N)$, where
    \begin{equation*}
        \Sigma_N = \bigcup_{j=1}^n (\xi_j + E[N]).
    \end{equation*}
    Then $\omega$ has no absolute saddle connections when $N \ge 2$. Indeed, if $\gamma$ were to be a saddle connection joining some $\xi = \xi_j + b \in \Sigma_N$ to itself, then developing to $\C$ yields a path
    \begin{equation*}
        \gamma(t) = \xi + (1-t)v_0 + tv_1,
    \end{equation*}
    for some $v_0,v_1 \in \Lambda$. But note that for $t = 1/N$, we have $\gamma(t) \in \xi_j + E[N]$, and so $\gamma$ has a marked point at time $t = 1/N$. Note that, assuming $\xi_1,...,\xi_n$ were chosen sufficiently generically, $\omega$ will have large orbit-closure.
\end{ex}

One also cannot replace the quadratic independence assumption of Lemma \ref{lemma:quadraticallyindependentholonomy} with the weaker assumption $\Hol_\omega = \Lambda_\Sigma$. To see this, note that in the proof of Lemma \ref{lemma:diophantinerelcomputation} we make nontrivial use of the quadratic independence assumption beyond the application of Lemma \ref{lemma:quadraticallyindependentholonomy}; namely, we use this assumption to remove exceptional subspaces. In fact, the conclusion of Lemma \ref{lemma:diophantinerelcomputation} does not follow from the weaker assumption $\Hol_\omega = \Lambda_\Sigma$, as illustrated by the following example.

\begin{ex}[Additional exceptional subspaces]\label{ex:extrasubspaces}
    Choose $\frac{p}{q} \in \Q \setminus \{0,1,-1\}$ in reduced form. For $\alpha,\beta \in \T^2$, let $\omega = (\T^2,\Sigma) \in \Hc(0,0^3)$ have marked points $\xi_1 = 0$, $\xi_2 = \alpha$, $\xi_3 = \beta$, $\xi_4 = \alpha + \frac{p}{q}\beta$. Assuming $\alpha,\beta$ are chosen sufficiently generically, one can show $\Hol_\omega = \Lambda_\Sigma$ and that the differences $\zeta_{ij}$ are distinct. However, the relation $q \zeta_{24} = p \zeta_{13}$ implies that many of the forms $Q_{ijkl} = Q_{\eta_{ijkl}}$ have exceptional subspaces, namely:
    \begin{align*}
        Q_{2413},Q_{4231} \rightsquigarrow \{(pz,qz)\}, \quad &\quad Q_{1324},Q_{3142} \rightsquigarrow \{(qz,pz)\}, \\
        Q_{2431},Q_{4213} \rightsquigarrow \{(pz,-qz)\}, \quad &\quad Q_{1342},Q_{3124} \rightsquigarrow \{(qz,-pz)\}.
    \end{align*}
    These subspaces contribute an additional $\sim \Delta_\Omega^{p,q} R^2$ to the count $N_{\Omega,A}(\omega,R)$, where
    \begin{equation*}
        \Delta_\Omega^{p,q} = 2\bigl(I_\Omega(p,q) + I_\Omega(q,p) + I_\Omega(p,-q) + I_\Omega(q,-p)\bigr), \quad \quad I_\Omega(a,b) = \int_{\R^2} \ind_\Omega(ax,bx)\dx.
    \end{equation*}
    For example, $\Delta_{\B_1^4}^{p,q} = \frac{8\pi}{p^2+q^2}$. Letting $\Ac^{p,q} \subseteq \Hc(0,0^3)$ be the subvariety consisting of tori whose marked points lie in the trapezoidal configuration above, we see $c_{\Omega,A}(\Ac^{p,q}) = F_{\Omega,1,4}(A) + \Delta_\Omega^{p,q}$. In particular, the map $\Ec \mapsto c_A(\Ec)$ will not be increasing under inclusion of subvarieties.
\end{ex}

\section{Branched covers}

In this section, we prove Theorem \ref{thm:1convergenceforbranchedcovers} for general subvarieties. To start, suppose that $\Ec = \Ec_{d,n}$ is a locus of degree-$d$ balanced covers of tori (possibly with marked points), branched over $n \ge 2$ points. Without loss of generality, we always assume one branch point is the origin. Each surface $\omega \in \Ec$ is equipped with a covering map $p \colon \omega \to (E,\dz)$, with $\omega = p^*\dz$. Note that $\Area(\omega) = d\Area(E)$.

\subsection{Diophantine surfaces and Theorem \ref{thm:1convergenceforbranchedcovers}}

There exists a \textit{forgetful map} $\Pi \colon \Ec \to \Hc(0,0^{n-1})$ sending $\omega$ to $(E,p(\Zeros(\omega)))$. By an argument analogous to Lemma \ref{lemma:combinatorialrefinement} (see also \cite[Lemma 2.2]{BarShu}), this map is of finite degree.\footnote{Here and throughout, we will ignore issues of orbifold points and coincident marked points.} Moreover, $\Pi$ commutes with the $\GL_2^+(\R)$-action. We record the following lemma.

\begin{lemma}\label{lemma:balancedholonomy}
    Suppose $\omega \in \Ec$. Then $\Hol_\omega = \Hol_{\Pi(\omega)}$.
\end{lemma}
\begin{proof}
    If $\gamma$ is a saddle connection on $\omega$, then its projection $\Bar{\gamma}$ to $\Pi(\omega)$ is a geodesic segment joining two marked points on $E$, with no marked points in its interior by the assumption on $\Ec$.

    Conversely, if $\gamma$ is a saddle connection on $\Pi(\omega)$ starting at a marked point $\xi_i$, then choose a preimage $\Tilde{\xi}_i$ for $\xi_i$ on $\omega$, and let $\Tilde{\gamma}$ be the trajectory starting at $\Tilde{\xi}_i$ of the same angle. Then at time $t_* = |\gamma|$, this saddle connection hits a preimage of the endpoint of $\gamma$, and since all preimages of branched points are marked on $\omega$ this defines a saddle connection so long as $\Tilde{\gamma}$ has no cone or marked points in its interior; but if it were to have a cone or marked point in its interior, then the corresponding geodesic $\gamma$ on $\Pi(\omega)$ would hit this marked point in its interior, a contradiction.
\end{proof}

\begin{defn}[Diophantine surfaces]\label{defn:diophantinesurf}
    Suppose $\omega \in \Ec$ has area $d$. We say that $\omega$ is \textit{(weakly) Diophantine} if $\Pi(\omega) \in \Hc^1(0,0^{n-1})$ is (weakly) Diophantine. If $\omega \in \Ec$ is of arbitrary area, we say it is (weakly) Diophantine if its area-$d$ rescaling is.
\end{defn}

Note that, if $\Ec$ is a locus of branched covers, then $\Ec$ contains a Diophantine surface if and only if $\Ec$ has full rel dimension $\rel \Ec = n-1$. In particular, this rules out any linear relations amongst the branch points. The following is an immediate consequence of Lemma \ref{lemma:balancedholonomy}.

\begin{cor}\label{cor:balancedholonomydiophantine}
    Suppose $\omega \in \Ec$ is a Diophantine surface with underlying torus $E = g\T^2$. Then $\Hol_\omega = g\Lambda_\Sigma$, where $\Sigma = p(\Zeros(\omega))$ and $\Lambda_\Sigma$ is as in \eqref{eq:quadraticallyindependentholonomy}.
\end{cor}

We also record a lemma, relating surfaces of general area to those of unit area.

\begin{lemma}\label{lemma:unitareaconversion}
    Suppose $\omega$ is a surface of area $v$, and $\omega' = v^{-1/2}\omega$ is its unit-area rescaling. Then
    \begin{equation*}
        c_{\Omega,A}(\omega') = v\cdot c_{\Omega,vA}(\omega).
    \end{equation*}
\end{lemma}
\begin{proof}
    This is a computation: we see that
    \begin{equation*}
        \frac{1}{R^2}N_{\Omega,A}(\omega',R) = \frac{1}{R^2}N_{\Omega,A}(v^{-1/2}\omega,R) = v \cdot \frac{1}{(v^{1/2}R)^2}N_{\Omega,vA}(\omega,v^{1/2}R),
    \end{equation*}
    and so taking a limit as $R \to \infty$ yields the desired result.
\end{proof}

The following result, combined with Example \ref{ex:kappalambdacomputation1}, generalizes Theorem \ref{thm:1convergenceforbranchedcovers}.

\begin{thm}\label{thm:5convergenceforbranchedcoversOmega}
    Suppose $\Ec$ is a degree-$d$ balanced locus of branched covers, branched over $n \ge 2$ points, and of full rel dimension. Then for every Diophantine $\omega \in \Ec^1$, every $\Omega \in \Starbodies$, and every $A \ge 0$,
    \begin{equation*}
        c_{\Omega,A}(\omega) = F_{\Omega,d,n}(A),
    \end{equation*}
    where $F_{\Omega,d,n}(A)$ is defined as in \eqref{eq:def-FOmegadnA}.
\end{thm}
\begin{proof}
    Suppose $\omega \in \Ec$ has underlying torus $E \in \Hc^1(0)$. By Corollaries \ref{cor:convergencefordiophantinetori} and \ref{cor:balancedholonomydiophantine}, we see that $c_{\Omega,A}(\omega) = F_{\Omega,1,n}(A)$. Since $\Area(\omega) = d$, Lemma \ref{lemma:unitareaconversion} completes the proof.
\end{proof}

\subsection{Subvarieties without additional markings}\label{subsec:nomarkedpts}

As noted in the introduction, the subvarieties of Theorem \ref{thm:1convergenceforbranchedcovers} require the balanced assumption, and hence not all loci of torus covers satisfy the necessary conditions. The goal of this section is to establish a combinatorial criterion for when a stratum $\Hc(\kappa)$ corresponding to a (positive) partition $\kappa$ admits a balanced locus.

Given partitions $\lambda,\mu$ of $N$, we say $\lambda$ \textit{refines} $\mu$ if the parts of $\mu$ can be further subdivided so as to form $\lambda$; a partition $\lambda$ is \textit{good} if $\lambda+1 = (\lambda_i+1)_i$ refines a rectangular partition $(d^n)$ for some $d,n \in \N$ with $d,n \ge 2$. Note every good partition is of length at least two.

\begin{prop}\label{prop:combinatorialclassification}
    Suppose $\kappa$ is a partition of $2g-2$. Then $\Hc(\kappa)$ admits a balanced subvariety $\Ec \subseteq \Hc(\kappa)$ if and only if $\kappa$ is good.
\end{prop}

Throughout, we let
\begin{equation*}
    \Pi_{1,n} = \generated{a,b,c_1,...,c_n}{aba^{-1}b^{-1}c_1\cdots c_n}
\end{equation*}
denote the abstract fundamental group of an $n$-times punctured torus.

\begin{lemma}\label{lemma:combinatorialrefinement}
    Suppose $\kappa$ is a partition of $2g-2$. Then $\Hc(\kappa)$ contains a balanced subvariety if and only if there exists $d,n \in \N$ with $2 \le d,n \le 2g-2$ and $dn = \dim_\C \Hc(\kappa) - 1$, and a representation $\rho \colon \Pi_{1,n} \to \Sym_d$, satisfying:
    \begin{enumerate}[label=(\roman*)]
        \item $\rho(\Pi_{1,n}) \le \Sym_d$ is a transitive subgroup;
        \item the monodromy $\rho(c_i)$ has no fixed points for each $i$;
        \item if $\rho(c_i) = \sigma_{i,1}\cdots \sigma_{i,s_i}$ is the cycle decomposition of $\rho(c_i)$, and $e_{ij} = \mathrm{length}(\sigma_{ij})$, then we have $(e_{ij} - 1)_{i,j} = \kappa$.
    \end{enumerate}
\end{lemma}
\begin{proof}
    Suppose $p \colon \omega \to (\T^2,\dz)$ is a Diophantine degree-$d$ branched cover of $\T^2$, branched over the set $\Sigma = \{\xi_1,...,\xi_n\} \subseteq \T^2$. Choose curves $a,b$ generating $\pi_1(\T^2)$, and for each $1 \le i \le n$, let $c_i$ be a small loop around the puncture $\xi_i$; this identifies $\pi_1(\T^2-\Sigma)$ with $\Pi_{1,n}$. Label the sheets of the cover by the integers $k \in \{1,...,d\}$. By the Galois correspondence (see e.g. \cite{Sza}), this data uniquely determines a \textit{monodromy representation} $\rho \colon \Pi_{1,n} \to \Sym_d$. Then $\omega$ being connected is equivalent to $\rho(\Pi_{1,n})$ being a transitive subgroup of $\Sym_d$. Moreover, for each $1 \le i \le n$, the cone angles of preimages of $\xi_i \in \Sigma$ under $p$ are encoded precisely in the cycle decomposition of $\rho(c_i)$; in particular, the preimage of $\xi_i$ on the $k$-th sheet of the cover is regular if and only if $\rho(c_i)(k) = k$, and it is straightforward to check that condition (iii) is equivalent to the statement that $\omega \in \Hc(\kappa)$. Therefore, a balanced locus $\Ec \subseteq \Hc(\kappa)$ exists if and only if there is \textit{some} $d \ge 2$ and $n \ge 2$ for which the conditions hold. It is easy to see that these conditions force $2 \le d,n \le 2g-2$; and moreover (iii) implies $dn = 2g-2+\mathrm{length}(\kappa) = \dim_\C \Hc(\kappa) - 1$.
\end{proof}

\begin{proof}[Proof of Proposition \ref{prop:combinatorialclassification}]
    By Lemma \ref{lemma:combinatorialrefinement}, it suffices to show a representation $\rho \colon \Pi_{1,n} \to \Sym_d$ satisfying (i)-(iii) exists if and only if $\kappa$ is good. It is clear that any partition $\kappa$ satisfying the conditions of Lemma \ref{lemma:combinatorialrefinement} is good for $(d^n)$, so we must demonstrate the converse.

    Suppose now that $\kappa$ is good for the partition $(d^n)$. For each $1 \le i \le n$, we have a subpartition $(\kappa_{j_1}+1,...,\kappa_{j_i}+1) \vdash d$; let $\rho(c_i) \in \Sym_d$ be some element with cycle type corresponding to this sub-partition. Since $\kappa_j \ge 1$, we clearly have (ii), and (iii) holds by definition. Moreover, so long as $n \ge 2$, it is straightforward to see that one can choose $\rho(c_1),\rho(c_2)$ to generate a transitive subgroup, so (i) is satisfied. It remains to show that $\rho$ can extend to a well-defined homomorphism on $\Pi_{1,n}$. This follows since
    \begin{equation*}
        \sgn(\rho(c_1)\cdots \rho(c_n)) = (-1)^{\sum_j [(\kappa_j+1)-1]} = (-1)^{2g-2} = 1,
    \end{equation*}
    and so $\rho(c_1)\cdots \rho(c_n) \in \Alt_d$; every element of $\Alt_d$ can be written as a $\Sym_d$-commutator by \cite[Theorem 1]{Ore}, so we can choose $\rho(a),\rho(b) \in \Sym_d$ for which $[\rho(a),\rho(b)]^{-1}=\rho(c_1)\cdots\rho(c_n)$; this produces the desired homomorphism.
\end{proof}

\begin{rmk}
    By Le Fils \cite[Theorem 1.3]{Fil}, one can choose $\Ec$ to lie in any \textit{component} of $\Hc(\kappa)$.
\end{rmk}

As it turns out, asymptotically almost every partition is good. To be precise, let $p(2g-2)$ denote the number of partitions of $2g-2$, and let $p_{\good}(2g-2)$ be the number of good partitions.

\begin{lemma}\label{lemma:aasgoodpartitions}
    We have
    \begin{equation*}
        \lim_{g \to \infty} \frac{p_{\good}(2g-2)}{p(2g-2)} =1.
    \end{equation*}
\end{lemma}
\begin{proof}
    Let $n = 2g-2$, let $\lambda = \kappa+1$, and $N = |\lambda|$. For a fixed prime $q$, we claim
    \begin{equation}\label{eq:partition1}
        \P\bigl(q \mid N,\, \text{$\lambda$ does not refine $((N/q)^q)$}\bigr) = o_{n \to \infty}(1).
    \end{equation}
    To see this, note that, for $n \gg_p 1$ there will be at least $p$ copies of both $2$ and $3$ in $\lambda$, with probability $1-o(1)$. Let $\Cs$ denote the multiset obtained by removing these copies of $2,3$ from $\lambda$; this implies
    \begin{equation}\label{eq:partition2}
        \sum_{\substack{a \in \Cs \\ a < x}} a \ge (q-1)x
    \end{equation}
    for every $x \ge 3$ appearing in $\Cs$, with probability $1-o(1)$. Indeed, for $x \le \epsilon \sqrt{n}$ this follows from the fact that there are at least $p\epsilon\sqrt{n}$ ones occurring in $\kappa$ with probability $1 - O_q(\epsilon)+o(1)$; on the other hand, Szalay-Tur\'an \cite[Theorem II]{SzaTur} implies that $\sum_{\epsilon \sqrt{n}/2 \le a \le \epsilon \sqrt{n}} a \gtrsim_\epsilon n$ with probability $1-o(1)$; finally $x \le n^{2/3}$ with probability $1-o(1)$, which completes the proof of \eqref{eq:partition2}.

    Now, take $q$ boxes $B_1,...,B_q$, and in non-increasing order add the elements of $\Cs$ to the lightest box. Up to reordering, we may assume that after this process has finished $B_1$ is the heaviest box. If the boxes are uneven, then the final element $x$ added to $B_1$ was a 2, since if $x \ge 3$ it would follow that $\sum_{a < x} a < (q-1)x$, contradicting \eqref{eq:partition2}. Thus, all boxes differ by weight at most 2.

    Now, add a 2 and a 3 to each box. If the boxes have uneven weight, then swap a 3 in the highest-weight box with a 2 in the lowest-weight one. Repeat this process until the boxes have weight differing by at most one; if $q \mid N$ then all boxes have equal weight. This establishes \eqref{eq:partition1}.

    Finally, note that the length of a random partition is asymptotically equidistributed in congruence classes (see e.g. \cite{BriFraGar}), so applying \eqref{eq:partition1} and sieving yields
    \begin{equation*}
        \limsup_{g \to \infty} \P\left( \kappa \vdash 2g-2 : \text{$\kappa$ is bad} \right) \le \limsup_{X \to \infty} \prod_{\substack{q \le X \\ \text{$q$ prime}}} \left(1 - \frac{1}{q} \right) = 0,
    \end{equation*}
    completing the proof.
\end{proof}

\begin{rmk}
    With more work one can show $\frac{p_{\good}(2g-2)}{p(2g-2)} = 1- O(1/\log g)$.
\end{rmk}

\begin{ex}[Rectangular partitions]\label{ex:factors}
    Suppose $2g-2 = k\ell$ for some $k,\ell \in \N$ with $\ell \ge 2$, and consider the partition $\kappa = (k^\ell)$. Set $d = k+1$, $n = \ell$, and let $g_1,...,g_n \in \Sym_d$ be $d$-cycles. Since $\kappa$ is rectangular, there exists a balanced locus $\Ec \subseteq \Hc(\kappa)$. In particular, balanced loci exist in all genera. The corresponding flat surfaces can be drawn as in the following figure.
    \begin{figure}[H]
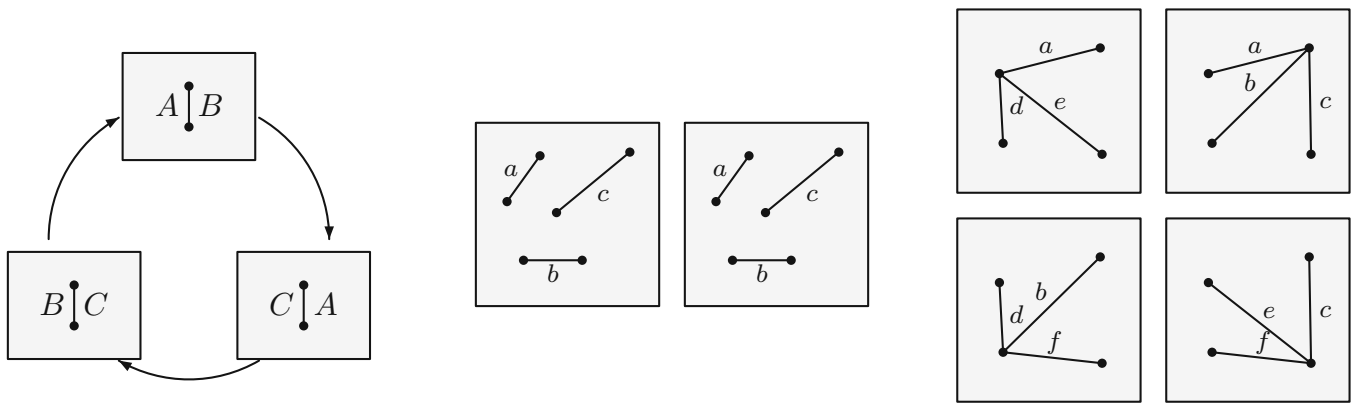

        \centering
        \BalancedCoversFigure[0.6]
        \label{fig:balancedcovers}
        \caption{Some branched covers corresponding to rectangular partitions}
    \end{figure}
    The case $g=2$, $k=1$, $\ell = 2$ produces $\Ec_4(1,1)$. More generally, the strata $\Hc(g-1,g-1)$ and $\Hc(1^{2g-2})$ are always good.
\end{ex}

\begin{ex}
    For a non-rectangular example, consider the partition $\kappa = (7,5,3,1)$. In this case, we can consider for example the representation
    \begin{equation*}
        \rho(c_1) = (1 \cdots 8)(9\, 10), \quad \quad \rho(c_2) = (1 \cdots 6)(7 \cdots 10).
    \end{equation*}
    On the other hand, the partition $\kappa = (7,5,4)$ is bad, since $8 + 6 + 5 = 19$ is prime.
\end{ex}

\begin{ex}[Different constants in the same stratum]
    The balanced covering $\Ec \subseteq \Hc(\kappa)$ produced in Proposition \ref{prop:combinatorialclassification} need not be unique, and in fact one can obtain two subvarieties $\Ec,\Ec' \subseteq \Hc(\kappa)$ in the same stratum with differing pair Siegel-Veech constants. For example, take $\kappa = (1^6)$. Let $\rho \colon \Pi_{1,2} \to \Sym_6$ be any representation extending the assignments
    \begin{equation*}
        \rho(c_1) = (12)(34)(56) \quad \mathrm{and} \quad \rho(c_2) = (13)(25)(46),
    \end{equation*}
    and let $\rho' \colon \Pi_{1,6} \to \Sym_2$ extend $\rho'(c_i) = (12)$ for all $i$. Then the balanced subvarieties $\Ec_\rho,\Ec_{\rho'}$ both lie in $\Hc(1^6)$, but it is clear that
    \begin{equation*}
        c_A(\Ec_\rho) = F_{6,2}(A) \ne F_{2,6}(A) = c_A(\Ec_{\rho'}).
    \end{equation*}
    One can moreover check that all of the individual terms of $F_{2,6}$ and $F_{6,2}$ disagree; thus, the pair Siegel-Veech constants associated to a subvariety do not seem to have an obvious relation to their ambient stratum.
\end{ex}

\begin{ex}[Same constants in different strata]
    Consider representations $\rho,\rho' \colon \Pi_{1,2} \to \Sym_6$ extending the assignments
    \begin{equation*}
        \rho(c_1) = (123)(456), \quad \rho(c_2) = (124)(356), \quad \quad \quad \rho'(c_1) = (12)(34)(56), \quad \rho'(c_2) = (123456).
    \end{equation*}
    Then the associated loci $\Ec_\rho \subseteq \Hc(2^4)$ and $\Ec_{\rho'} \subseteq \Hc(5,1^3)$ have the same pair Siegel-Veech constants in the same genus $g = 5$, but lie in different strata.
\end{ex}

\subsection{Siegel-Veech transforms}\label{subsec:SVtransforms}

Suppose $\Ec$ is an affine subvariety. By Veech \cite{Vee1}, there exists a constant $c_{\SV}(\Ec) > 0$ so that for all $f \in C_c(\R^2)$,
\begin{equation*}
    \int_{\Ec^1} \sum_{z \in \Hol_\omega} f(z)\dmu_\Ec(\omega) = c_{\SV}(\Ec) \int_{\R^2} f\dx.
\end{equation*}
The asymptotic \eqref{eq:largeAlimit} yields the following immediate corollary.

\begin{cor}\label{cor:cAEsquareSV}
    If $\Ec$ is as in Theorem \ref{thm:1convergenceforbranchedcovers}, then
    \begin{equation*}
        c_{\Omega,0}(\Ec) = \kappa_\Omega \cdot  c_{\SV}(\Ec) \quad \mathrm{and} \quad \lim_{A \to \infty} \frac{c_{\Omega,A}(\Ec)}{A} = 2\lambda_\Omega \cdot c_{\SV}(\Ec)^2.
    \end{equation*}
\end{cor}
\begin{proof}
    If $\Ec$ has full rel dimension, then a generic surface $\omega \in \Ec$ does not have any non-homologous parallel vectors. It is then straightforward to see $N_{0}(\omega,R) = 2|\Hol_\omega \cap \B_{R/\sqrt{2}}^2|$, after which the first equality follows from \cite{Vee1}. For the second equality, apply \eqref{eq:largeAlimit}.
\end{proof}

Using the Riesz representation theorem and the decomposition of $\C \times \C$ as an $\SL_2(\R)$-representation, Athreya-Cheung-Masur \cite{AthCheMas} showed there are measures $\eta_\Ec$ on $\P_\R^1$ and $\nu_\Ec$ on $\R^\times$ so that
\begin{equation*}
    \int_{\Ec^1} \sum_{z,w \in \Hol_\omega} f(z,w) \dmu_\Ec = \int_{\P_\R^1} \int_{\R^2} f(x,sx)\dx\deta_\Ec(s) + \int_{\R^\times} \int_{\Sigma_t} f(v) \dsigma_t(v)\dnu_\Ec(t)
\end{equation*}
for all $f \in C_c(\R^2 \times \R^2)$. For $\Ec$ as in Theorem \ref{thm:1convergenceforbranchedcovers}, we are able to compute explicitly the measures $\eta_\Ec$ and $\nu_\Ec$, yielding an infinite family of (non-Veech) subvarieties for which \cite[Question 4.4]{AthCheMas} has an explicit answer.

\begin{prop}
    Suppose $\Ec$ is as in Theorem \ref{thm:1convergenceforbranchedcovers}. Then
    \begin{equation*}
        \eta_\Ec = \bigl(6d/\pi^2 + dn(n-1)\bigr)\bigl[ \delta_{-1} + \delta_{+1} \bigr],
    \end{equation*}
    and
    \begin{equation*}
        \nu_\Ec = \bigl(12d^2n(n-1)/\pi^2 + d^2n^2(n-1)^2 \bigr)|t|\dt + \frac{6}{\pi^2}\sum_{\substack{t \in \Z \\ t \ne 0}} \phi(t)\delta_{t/d}.
    \end{equation*}
\end{prop}

\begin{cor}
    For the eigenform locus $\Ec_4(1,1)$,
    \begin{align*}
        \eta_{\Ec_4(1,1)} &= \left( 12/\pi^2 + 4 \right)\bigl[\delta_{-1}+\delta_{+1}\bigr] \\
        \nu_{\Ec_4(1,1)} &= \left( 96/\pi^2 + 16 \right)|t|\dt + \frac{6}{\pi^2} \sum_{\substack{t \in \Z \\ t \ne 0}} \phi(t)\delta_{t/2}.
    \end{align*}
\end{cor}

\subsection{Rectangular billiards}\label{subsec:billiards}

\begin{figure}[H]
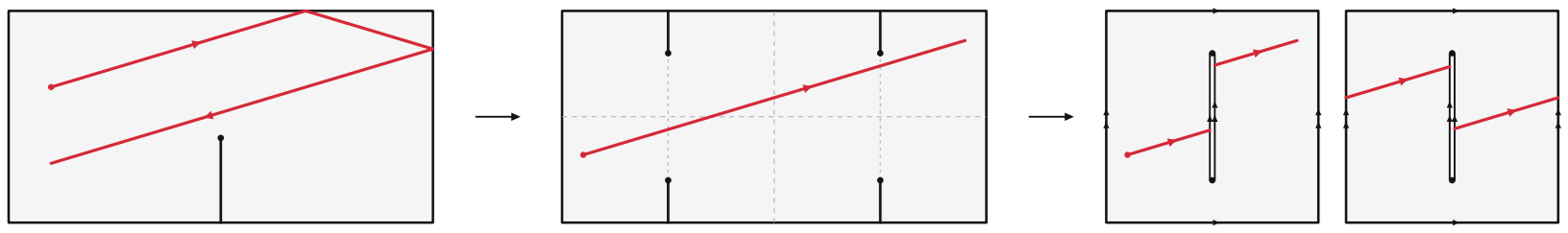

    \centering
    \BilliardUnfoldingFigure[0.75]
    \caption{Unfolding rectangular billiards with a barrier of length $\frac{1-\alpha}{2}$, with $\alpha = \sqrt{2}-1$}
    \label{fig:unfolding}
\end{figure}

Suppose $\alpha \in (0,1)$, and let $P_\alpha$ be a rectangle of dimensions $1 \times \frac{1}{2}$, with a vertical barrier of length $\frac{1-\alpha}{2}$ emanating from the middle of a horizontal side. By a well-known unfolding procedure (see e.g. \cite{EskMasSch}, \cite[Chapter 4]{AthMas}), the billiard flow on $P_\alpha$ with elastic collisions is equivalent to the straight-line flow on the slit torus $\omega = \T^2 \#_\alpha \T^2$ obtained by taking two square tori and gluing them along a vertical slit from $0$ to $\alpha$. See figure \ref{fig:unfolding}.

Note that $\omega$ will not be Diophantine in the sense of Definition \ref{defn:diophantinesurf}, since its holonomy set misses (essentially) the entire set of $\Z$-lines
\begin{equation*}
    V_\pm = \bigl\{(0,k\pm \alpha) : k \in \Z \setminus \{0\} \bigr\}
\end{equation*}
from the set of saddle connections, which has the potential to affect pair counting by a meaningful factor of $\asymp R^2$. For example, the contribution of $V_\pm$ is essentially given by the count
\begin{equation*}
    \bigl|\left\{ x \in \Z^4 : \norm{x} \le R,\ x_1 = 0,\ a < Q_{\pm\eta}(x) < b \right\}\bigr|
\end{equation*}
for $\eta = (0,\alpha,0,\alpha)$. Expanding this out, we are reduced to counting integer solutions to the system of inequalities
\begin{equation*}
    |y(z \pm\alpha)| < A \quad \mathrm{and} \quad
    y^2 + z^2 + w^2 \le R^2
\end{equation*}
for some fixed $A \ge 1$ and $R \gg_A 1$. Setting $y = 0$ yields a contribution of $\sim 2\pi R^2$ to the counting function from these two lines. 

On the other hand, the lines $V_\pm$ correspond to exceptional subspaces for $Q_{\pm \eta}$, and so Margulis-Mohammadi implies that all counting asymptotics are valid only once appropriate shifts of integer vectors lying along $V_\pm$ have been removed from the counting function; in particular, the set $\Lambda_{\Sigma}$ defined in \eqref{eq:quadraticallyindependentholonomy} is no longer the correct object of study, but rather $\Lambda_{\Sigma} \setminus V_\pm$, which is indeed $\Hol_\omega$. Therefore, the proof of Theorem \ref{thm:1convergenceforbranchedcovers} easily extends to show the following.

\begin{prop}
    Suppose $\omega = \T^2 \#_\alpha \T^2$, with $\alpha$ Diophantine. Then for all $A \ge 0$ we have
    \begin{equation*}
        c_{A}(\omega) = \frac{6}{\pi} + 2\pi + 12\Phi_3(A) + (8\pi^2 + 48)A.
    \end{equation*}
\end{prop}

We also illustrate how to derive the \textit{multiset} constant $\Tilde{c}_{\Omega,A}$ obtained by viewing $\Hol_\omega$ as a set with multiplicity. For $\omega = \T^2 \#_\alpha \T^2$, it is straightforward to see that the weighted decomposition of \eqref{eq:quadraticallyindependentholonomy} becomes
\begin{equation*}
    \widetilde{\Hol}_\omega = 4\{\Z_{\prim}^2 \setminus V_{\pm} \}+2\{(\Z^2+(0,\alpha)) \setminus V_\pm \} + 2\{(\Z^2-(0,\alpha)) \setminus V_\pm \}.
\end{equation*}
Therefore, we obtain analogously to \eqref{eq:cOmegaAdecomp} the equality
\begin{equation*}
    \Tilde{c}_{\Omega,A}(\omega) = 16c_{\Omega,A}^{\abs}(\omega) + 4\sum_{\substack{1 \le i,j,k,l \le 2 \\ i \ne j, k \ne l}} c_{\Omega,A}^{\rel,ijkl}(\omega) + 8\sum_{\substack{1 \le i,j \le 2 \\ i \ne j}} c_{\Omega,A}^{\mix,ij}(\omega),
\end{equation*}
where $c_{\Omega,A}^*(\omega)$ is defined as in Section \ref{subsec:diophantine}.

\begin{prop}\label{prop:multiset}
    Suppose $\omega = \T^2 \#_\alpha \T^2$, with $\alpha$ Diophantine. Then for all $\Omega \in \Starbodies$ and $A \ge 0$,
    \begin{equation*}
        \Tilde{c}_{\Omega,A}(\omega) = \frac{16\kappa_\Omega}{\zeta(2)} + 8\kappa_\Omega + \frac{32\lambda_\Omega}{\zeta(2)}\Phi_3(A) + \lambda_\Omega\biggl(32 + \frac{64}{\zeta(2)} \biggr)A.
    \end{equation*}
\end{prop}

Note the \textit{weighted} Siegel-Veech constant $\Tilde{c}_{\SV}(\Ec_4(1,1))$ is given by
\begin{equation*}
    \Tilde{c}_{\SV}(\Ec_4(1,1)) = 8(1 + \zeta(2)^{-1}),
\end{equation*}
and so once again we have
\begin{equation}\label{eq:svsquared3}
    \lim_{A \to \infty} \frac{\Tilde{c}_{\Omega,A}(\Ec_4(1,1))}{A} =  2\lambda_\Omega \cdot \Tilde{c}_{\SV}(\Ec_4(1,1))^2.
\end{equation}

\section{Veech surfaces}\label{sec:veech}

In this section, we demonstrate how the proof of Lemma \ref{lemma:diophantineabscomputation} can be modified to prove Theorem \ref{thm:3convergenceforveech}, using results of Burrin-Fairchild \cite{BurFai} on Fuschsian groups. Given a Veech curve $\Vc \subseteq \Hc$, we define its set of \textit{virtual triangles} $\VT(\Vc)$ by
\begin{equation*}
    \VT(\Vc) = \left\{ |z \wedge w| : z,w \in \Hol_\omega \right\},
\end{equation*}
where $\omega \in \Vc^1$ is arbitrary. It is easy to see that $\VT(\Vc)$ is independent of the choice of $\omega$. By Smillie-Weiss \cite{SmiWei}, this is a discrete subset of $[0,\infty)$. Our goal is to prove the following generalization of Theorem \ref{thm:3convergenceforveech}.

\begin{thm}\label{thm:6convergenceforveechOmega}
    Suppose $\Vc$ is a Veech curve. Then there exists some $\kappa_{\Omega,\Vc} > 0$ and some integer-valued weights $\phi_\Vc \colon \VT(\Vc) \to \N$ so that for every $\omega \in \Vc^1$,
    \begin{equation*}
        c_{\Omega,A}(\omega) = \frac{\kappa_{\Omega,\Vc}}{\Vol(\Vc^1)} + \frac{2\lambda_\Omega}{\Vol(\Vc^1)}\sum_{t \in \VT(\Vc) \cap (0,A]} \frac{\phi_\Vc(t)}{t}.
    \end{equation*}
\end{thm}

To begin the proof, fix $\omega \in \Vc^1$ and let $\Gamma = \Stab_{\SL_2(\R)}(\omega)$. For the sake of simplicity, we will assume $-\Id \in \Gamma$; for the general case, take $\Tilde{\Gamma} = \gen{\Gamma,-\Id}$, run the proof for $\Tilde{\Gamma}$, and then use the relation $\Covol(\Tilde{\Gamma}) = \frac{1}{2}\Covol(\Gamma)$.

By Veech \cite{Vee}, the set $\Hol_\omega$ is a finite union of discrete $\Gamma$-orbits corresponding to the cusps $\Cusps(\Gamma) = \del(\H^2/\Gamma)$ of $\Gamma$. For each cusp $\afrak \in \Cusps(\Gamma)$, let $u_\afrak \in \S^1$ be a choice of direction for the cusp, and $h_\afrak$ its width. Let $z_{\afrak,1},...,z_{\afrak,r_\afrak}$ denote the vectors in $\Hol_\omega$ in direction $u_\afrak$. Then $z_{\afrak,j} = \ell_{\afrak,j}u_\afrak$ for some $\ell_{\afrak,1} > \cdots > \ell_{\afrak,r_\afrak} > 0$, so
\begin{equation}\label{eq:Veechholonomydecomp}
    \Hol_\omega = \bigcup_{\afrak \in \Cusps(\Gamma)}\bigcup_{j=1}^{r_\afrak} \ell_{\afrak,j}\Gamma u_\afrak.
\end{equation}
For each $t \in \VT(\Vc)$, each $\afrak,\bfrak \in \Cusps(\Gamma)$, and each $1 \le i \le r_\afrak$, $1 \le j \le r_\bfrak$, we set
\begin{equation*}
    \Bar{N}_{\Omega,t}^{\afrak,\bfrak;i,j}(\omega,R) = \bigl| \bigl\{ (z,w) \in \ell_{\afrak,i}\Gamma u_\afrak \times \ell_{\bfrak,j}\Gamma u_\bfrak : (z,w) \in R\Omega,\, |z \wedge w| = t \bigr\} \bigr|.
\end{equation*}
By \eqref{eq:Veechholonomydecomp}, it follows that
\begin{equation}\label{eq:VeechNAdecomp}
    N_{\Omega,A}(\omega,R) = \sum_{t \in \VT(\Vc) \cap [0,A]} \sum_{\afrak,\bfrak \in \Cusps(\Gamma)} \sum_{\substack{1 \le i \le r_\afrak \\ 1 \le j \le r_\bfrak}} \Bar{N}_{\Omega,t}^{\afrak,\bfrak;i,j}(\omega,R).
\end{equation}

\begin{lemma}\label{lemma:veechc0}
    For each $\Omega \in \Starbodies$, we have
    \begin{equation*}
        c_{\Omega,0}(\omega) = \frac{\kappa_{\Omega,\omega}}{\Vol(\Vc^1)}
    \end{equation*}
    for some constant $\kappa_{\Omega,\omega} > 0$.
\end{lemma}
\begin{proof}
    The only contributions to the count $N_{\Omega,0}(\omega,R)$ come from pairs corresponding to the same cusp, and so by \eqref{eq:VeechNAdecomp} we have
    \begin{equation*}
        N_{\Omega,0}(\omega,R) = \sum_{\afrak \in \Cusps(\Gamma)} \sum_{1 \le i,j \le r_\afrak} \Bar{N}_{\Omega,0}^{\afrak,\afrak;i,j}(\omega,R).
    \end{equation*}
    Let us introduce the notation $T_{a,b} = \diag(a^{-1},b^{-1})$, viewed as acting on the right of $\C \times \C$ by multiplication. Then
    \begin{equation*}
        \Bar{N}_{\Omega,0}^{\afrak,\afrak;i,j}(\omega,R) = \bigl| \bigl\{ (\gamma u_\afrak,\pm \gamma u_\afrak) : \gamma \in \Gamma,\, (\gamma u_\afrak,\pm \gamma u_\afrak) \in R\Omega. T_{\ell_{\afrak,i},\ell_{\afrak,j}} \bigr\} \bigr|.
    \end{equation*}
    Once again applying \cite{EskMcM} with appropriate normalizations, we see that
    \begin{equation*}
        \Bar{N}_{\Omega,0}^{\afrak,\afrak;i,j}(\omega,R) \sim \frac{h_\afrak}{\Covol(\Gamma)}R^2\int_{\R^2} \biggl[\ind_\Omega(\ell_{\afrak,i}x, \ell_{\afrak,j}x) + \ind_\Omega(\ell_{\afrak,i}x,-\ell_{\afrak,j}x)\biggr]\dx.
    \end{equation*}
    Taking a limit and summing over $\afrak,i,j$, we obtain $c_{\Omega,0}(\omega) = \frac{\kappa_{\Omega,\Vc}}{\Vol(\Vc^1)}$, where
    \begin{equation}\label{eq:kappaveech}
        \kappa_{\Omega,\Vc} = \sum_{\afrak \in \Cusps(\Gamma)} h_\afrak \sum_{1 \le i,j \le r_\afrak} \int_{\R^2} \biggl[\ind_\Omega(\ell_{\afrak,i}x, \ell_{\afrak,j}x) + \ind_\Omega(\ell_{\afrak,i}x,-\ell_{\afrak,j}x)\biggr]\dx.
    \end{equation}
\end{proof}

\begin{lemma}\label{lemma:veechct}
    For each $t \in \VT(\Vc)$ nonzero, each $\Omega \in \Starbodies$, each $\afrak,\bfrak \in \Cusps(\Gamma)$, and each $1 \le i \le r_\afrak$, $1 \le j \le r_\bfrak$, there exists some $\phi_{\omega,\afrak,\bfrak,i,j}(t) \in \Z_{\ge 0}$ for which
    \begin{equation*}
        \lim_{R \to \infty} R^{-2}\Bar{N}_{\Omega,t}^{\afrak,\bfrak;i,j}(\omega,R) = \frac{2\lambda_\Omega \phi_{\omega,\afrak,\bfrak,i,j}(t)}{t\Covol(\Gamma)}.
    \end{equation*}
    For some choice of $\afrak,\bfrak,i,j$, this integer is positive.
\end{lemma}
\begin{proof}
    Let $\phi_{\omega,\afrak,\bfrak,i,j}(t) \in \Z_{\ge 0}$ denote the number of $\Gamma$-orbits of pairs of holonomy vectors $(z,w) \in \ell_{\afrak,i}\Gamma u_\afrak \times \ell_{\bfrak,j}\Gamma u_\bfrak$ with $z \wedge w = t$; by \cite[Lemma 7.2]{BurFai}, this number is finite.\footnote{\cite{BurFai} obtain asymptotics for these counts.} Each of these orbits has a representative of size $O_{\Gamma,t}(1)$, and so for $R \gg_{\Gamma,t} 1$ we may argue exactly as in Lemma \ref{lemma:diophantineabscomputation} to obtain
    \begin{equation*}
        \Bar{N}_{\Omega,t}^{\afrak,\bfrak;i,j}(\omega,R) \sim \frac{2\lambda_\Omega\phi_{\omega,\afrak,\bfrak,i,j}(t)}{t\Covol(\Gamma)}R^2,
    \end{equation*}
    as desired.
\end{proof}

\begin{proof}[Proof of Theorem \ref{thm:6convergenceforveechOmega}]
    Combining Lemmas \ref{lemma:veechc0} and \ref{lemma:veechct}, we see that for $\omega \in \Vc^1$ fixed we have
    \begin{equation*}
        c_{\Omega,A}(\omega) = \frac{\kappa_{\Omega,\omega}}{\Vol(\Vc^1)} + \frac{2\lambda_\Omega}{\Vol(\Vc^1)}\sum_{t \in \VT(\Vc) \cap (0,A]} \frac{\phi_\omega(t)}{t},
    \end{equation*}
    where
    \begin{equation}\label{eq:Wveech}
        \phi_{\omega}(t) = \sum_{\afrak,\bfrak \in \Cusps(\Gamma)} \sum_{\substack{1 \le i \le r_\afrak \\ 1 \le j \le r_\bfrak}} \phi_{\omega,\afrak,\bfrak,i,j}(t).
    \end{equation}
    Finally, note that both $\kappa_{\Omega,\omega}$ and $\phi_\omega$ are $\SL_2(\R)$-invariant.
\end{proof}

As an immediate corollary, we see that the constants $c_{\Omega,A}(\Vc)$ are ``arithmetic,'' in the following sense. Let $\kay_\Vc = \Q(\tr \Gamma)$ denote the \textit{trace field} of $\Vc$.

\begin{cor}
    For every $\Omega \in \Starbodies$ and every $A \ge 0$,
    \begin{equation*}
        c_{\Omega,A}(\Vc) \in \kay_\Vc\bigl( \Vol(\Vc^1), \kappa_{\Omega,\Vc},\lambda_\Omega \bigr).
    \end{equation*}
\end{cor}
\begin{proof}
    By Kenyon-Smillie \cite[Theorems 28-30]{KenSmi} or McMullen \cite[Theorems 9.4-9.5, Corollary 9.6]{McM}, for any $\omega \in \Vc$ there exists $g \in \GL_2^+(\R)$ so that $\hol_{g\omega}(H_1(X,\Sigma;\Z)) \subseteq \kay_\Vc + i\kay_\Vc$; since the area of a surface can be written as a quadratic form on the absolute periods, it follows that $\Area(g\omega) \in \kay_\Vc$, and so we can apply an element $h \in \GL_2^+(\kay_\Vc)$ so that $\omega' = hg\omega \in \Vc^1$ and $\hol_{\omega'}(H_1(X,\Sigma;\Z)) \subseteq \kay_\Vc + i\kay_\Vc$. Then $\VT(\Vc) = \VT(\omega') \subseteq \kay_\Vc$, which combined with Theorem \ref{thm:6convergenceforveechOmega} completes the proof.
\end{proof}

\begin{rmk}
    The bound $\limsup_{R \to \infty} \frac{1}{R^2}N_A(\omega,R) \le C_\omega A + D_\omega$ for some constants $C_\omega,D_\omega > 0$ depending on $\omega$ was previously shown in \cite[Theorem 1.6]{BurFai}.
\end{rmk}

\begin{rmk}
    An alternative proof of Theorem \ref{thm:3convergenceforveech}, which is ``more geometric'' but does not recover the exact constants, can be given as follows. The zero-determinant contribution reduces to Eskin-Masur \cite{EskMas}. For $t \ne 0$ one can define a restricted pair Siegel-Veech transform $\Theta_t$ on $C_c(\Sigma_t)$, summing over pairs of holonomy vectors with determinant $t$. Since $\Sigma_t$ is a principal $\SL_2(\R)$-torsor, the map $f \mapsto \int_{\Vc^1} \Theta_t f \dmu_\Vc$ is a multiple of the Lebesgue measure $\sigma_t$. A renormalization argument along the level sets $\Sigma_t$ using a variant of the trapezoid lemma (\cite[Lemma 3.6]{EskMarMoz}) for pairs, alongside Ratner's theorems and Smillie-Weiss \cite{SmiWei}, completes the proof.
\end{rmk}

\subsection{Hecke triangle groups}\label{subsec:heckegroups}

In this section, we specialize to the case of Veech curves arising from regular polygons. These surfaces were first studied by Veech \cite{Vee}.

To start, fix some $p \in \N$, let $q = 2p+1$ be odd, and set
\begin{equation*}
    \theta = 2\pi/q, \quad \quad \zeta = e^{i\theta}, \quad \quad h = 2\cot(\theta/2), \quad \quad Z_q = \sum_{1 \le i,j \le p} \frac{1}{\sin^2(i\theta)+\sin^2(j\theta)}.
\end{equation*}
Note that $Z_q \in \Q$, since it is $\Gal(\Q(\zeta)/\Q)$-invariant.

Let $\Delta \subseteq \C$ be the triangle with vertices $0,1,\zeta$, and let $P = \bigcup_{j=0}^{q-1}\zeta^j\Delta$. By gluing $P$ along a reflected copy $Q$, we obtain a translation surface $X_q$ of genus $p$, lying in $\Hc(2p-2)$.

\begin{figure}[H]
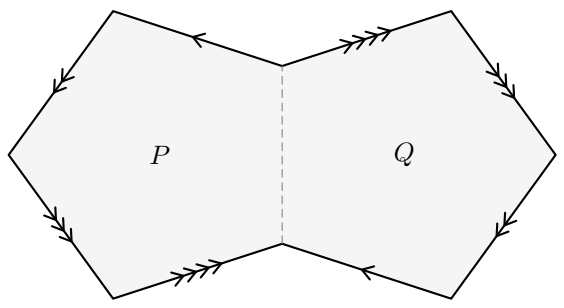

    \centering
    \DoublePentagonFigure[0.5]
    \label{fig:pentagon}
    \caption{The surface $X_5$}
\end{figure}

\begin{thm}[\cite{Vee}, Theorem 5.8]
    The surface $X_q$ is Veech, and its Veech group $\Gamma_q$ is generated by the matrices
    \begin{equation*}
        \begin{bmatrix}
            \cos \theta & -\sin \theta \\
            \sin \theta & \cos \theta
        \end{bmatrix} \quad \mathrm{and} \quad \begin{bmatrix}
            1 & 0 \\
            h & 1
        \end{bmatrix}.
    \end{equation*}
    The quotient $\Gamma_q/\{\pm \Id\}$ is isomorphic to the $(2,q,\infty)$-triangle group. In particular, $\Gamma_q$ has a single cusp $\ofrak$, with $h_\ofrak = h$.
\end{thm}

Under our normalization, the cusp $\ofrak$ corresponds to $0 \in \del \H^2$. Note $\Area(X_q) = q\sin(\theta)$; let $\omega_q = \Area(X_q)^{-1/2}X_q$ be its unit-area rescaling. We write $\Vc_q = \GL_2^+(\R)X_q$ for the orbit of $X_q$ in $\Hc(2p-2)$, with volume $\Vol(\Vc_q^1) = \frac{\pi^2(q-2)}{2q}$.

\begin{lemma}\label{lemma:2qinftyholonomy}
    The set $\Hol_{\omega_q}$ decomposes as
    \begin{equation*}
        \Hol_{\omega_q} = \bigcup_{j=1}^p \ell_j \Gamma_q \ebf_2,
    \end{equation*}
    where $\ell_j = 2(q\sin(\theta))^{-1/2} \Im(\zeta^j) = 2(q\sin(\theta))^{-1/2}\sin(j\theta)$, and $\ebf_2 \cong i$ is the vertical unit vector.
\end{lemma}
\begin{proof}
    It is clear from the polygonal construction of $X_q$ that $\{\ell_j \ebf_2\}_{j=1}^p$ constitutes the full set of \textit{vertical} holonomy vectors on $X_q$, after which the decomposition follows from the fact that $\Vc_q$ has a single cusp.
\end{proof}

We would like to relate the set $\VT(\Vc_q)$ to the \textit{$q$-geometric totient function}, as defined in Fairchild \cite[Section 2.4]{Fai}. To do so, we introduce the \textit{standard Hecke triangle group $H_q$}, given by
\begin{equation*}
    H_q = \left\langle \begin{bmatrix}
        0 & -1 \\
        1 & 0
    \end{bmatrix}, \begin{bmatrix}
        1 & 0 \\
        2\cos(\theta/2) & 1
    \end{bmatrix} \right\rangle
\end{equation*}
By Arnoux-Hubert \cite{ArnHub} (see also Arnoux-Schmidt \cite{ArnSch}), the matrix
\begin{equation*}
    U = \begin{bmatrix}
        1 & 0 \\
        \cos(\theta/2) & \sin(\theta/2)
    \end{bmatrix}
\end{equation*}
satisfies $U\Gamma_qU^{-1} = H_q$. Let $\Tc_q = \{|z \wedge w| : z,w \in H_q\ebf_2\}$. Following \cite{Fai},\footnote{\cite{Fai} defines the totient function with respect to a horizontally based orbit, but since $H_q$ contains a rotation by $\pi/2$ these definitions coincide.} we define the \textit{$q$-geometric totient function} $\phi_q$ by
\begin{equation*}
    \phi_q(s) = \left|\left\{ 1 \le a \le 2\cos(\pi/q)s : s+ia \in H_q\ebf_2 \right\} \right|.
\end{equation*}

\begin{lemma}[\cite{Fai}, Lemma 3.6]\label{lemma:Fai36}
    For each $s \in \Tc_q$ nonzero, the set of pairs $(z,w) \in H_q\ebf_2 \times H_q\ebf_2$ with $z \wedge w = s$ is a union of $\phi_q(s)$ orbits of $H_q$.
\end{lemma}

Lemma \ref{lemma:2qinftyholonomy} implies the following.

\begin{lemma}\label{lemma:VTq}
    The set of virtual triangle areas of $\Vc_q$ is given by
    \begin{equation*}
        \VT(\Vc_q) = \sin(\theta/2)\bigcup_{1 \le i,j \le p} \ell_i\ell_j\Tc_q.
    \end{equation*}
\end{lemma}

Let us define the \textit{$q$-geometric summatory function} $\Phi_q$ by
\begin{equation*}
    \Phi_q(T) = \sum_{\substack{t \in \Tc_q \\ 0 < t \le T}} \frac{\phi_q(t)}{t}.
\end{equation*}
Note that $\Phi_q$ with $q=3$ agrees with our earlier definition $\Phi_3(T) = \sum_{t=1}^{\floor{T}} \frac{\phi(t)}{t}$.

\begin{prop}\label{prop:cAqveech}
    We have
    \begin{equation*}
        c_{A}(\Vc_q) = \frac{4q^2\cos^2(\theta/2)Z_q}{\pi(q-2)} + \frac{2q^2\cos(\theta/2)}{q-2} \sum_{1 \le i,j \le p} \frac{1}{\sin(i\theta)\sin(j\theta)} \Phi_q\biggl( \frac{A}{\sin(\theta/2)\ell_i\ell_j} \biggr).
    \end{equation*}
\end{prop}
\begin{proof}
    This is an application of Theorem \ref{thm:6convergenceforveechOmega}. To start, in the case $\Omega = \B_1^4$ one can compute $\kappa_{\Vc_q}$ as in \eqref{eq:kappaveech}:
    \begin{equation*}
        \kappa_{\Vc_q} = 2\cot(\theta/2) \sum_{1 \le i,j \le p} \frac{2\pi}{\ell_i^2 + \ell_j^2} = \pi q\sin(\theta)\cot(\theta/2)Z_q.
    \end{equation*}
    In view of Lemma \ref{lemma:VTq}, we have
    \begin{equation}\label{eq:qcurvedecomp}
        \sum_{1 \le i,j \le p} \sum_{t \in (0,A] \cap \sin(\theta/2)\ell_i\ell_j\Tc_q} \Bar{N}_{t}^{\ofrak,\ofrak;i,j}(\omega,R) = \sum_{1 \le i,j \le p} \sum_{s \in (0,\frac{A}{\sin(\theta/2)\ell_i\ell_j}] \cap \Tc_q} \Bar{N}_{\sin(\theta/2)\ell_i\ell_js}^{\ofrak,\ofrak;i,j}(\omega,R).
    \end{equation}
    By combining Lemma \ref{lemma:veechct} with Lemma \ref{lemma:Fai36},
    \begin{equation*}
        \Bar{N}_{\sin(\theta/2)\ell_i\ell_js}^{\ofrak,\ofrak;i,j}(\omega,R) \sim \frac{2\lambda_{\B_1^4}\phi_q(s)}{\sin(\theta/2)\ell_i\ell_j s \Covol(\Gamma_q)}R^2.
    \end{equation*}
    Thus, we have
    \begin{align*}
        \eqref{eq:qcurvedecomp} &\sim R^2\frac{2\lambda_{\B_1^4}}{\sin(\theta/2)\Covol(\Gamma_q)} \sum_{1 \le i,j \le p} \sum_{s \in (0,\frac{A}{\sin(\theta/2)\ell_i\ell_j}] \cap \Tc_q} \frac{\phi_q(s)}{\ell_i\ell_j s} \\
        &\sim R^2 \frac{q^2\sin(\theta)}{\sin(\theta/2)(q-2)} \sum_{1 \le i,j \le p} \frac{1}{\sin(i\theta)\sin(j\theta)} \sum_{s \in (0,\frac{A}{\sin(\theta/2)\ell_i\ell_j}] \cap \Tc_q} \frac{\phi_q(s)}{s},
    \end{align*}
    completing the proof.
\end{proof}

\begin{cor}
    As $A \to \infty$, we have
    \begin{equation*}
        c_A(\Vc_q) \sim \frac{2q^4(q^2-1)^2\cos^4(\pi/q)}{9\pi^2(q-2)^2}A.
    \end{equation*}
\end{cor}
\begin{proof}
    By \cite[Theorem 7.4]{BurFai} (which itself is a corollary of \cite[Chapter 11, Theorem 4]{Goo}), alongside a straightforward integration by parts argument, we see that
    \begin{equation*}
        \Phi_q(T) \sim \frac{4\cos^2(\theta/2)}{\Vol(\Vc_q^1)}T.
    \end{equation*}
    The large-$A$ asymptotic then follows from Proposition \ref{prop:cAqveech}:
    \begin{align*}
        \frac{c_A(\Vc_q)}{A} &\sim \frac{1}{A}\frac{2\pi^2}{\sin(\theta/2)\Vol(\Vc_q^1)}\sum_{1 \le i,j \le p} \frac{1}{\ell_i\ell_j} \Phi_q\biggl( \frac{A}{\sin(\theta/2)\ell_i\ell_j} \biggr) \\
        &\sim \frac{2\pi^2(2\cos(\theta/2))^2}{\sin^2(\theta/2)\Vol(\Vc_q^1)^2} \sum_{1 \le i,j \le p} \frac{1}{\ell_i^2\ell_j^2} \\
        &=  \frac{8\pi^2\cos^2(\theta/2)q^2\sin^2(\theta)}{16\sin^2(\theta/2)\Vol(\Vc_q^1)^2} \sum_{1 \le i,j \le p} \frac{1}{\sin^2(i\theta)\sin^2(j\theta)} \\
        &= \frac{32\pi^2\cos^2(\theta/2)\sin^2(\theta)q^4}{16\sin^2(\theta/2)\pi^4(q-2)^2}\left( \sum_{j=1}^p \csc^2(j\theta) \right)^2 \\
        &= \frac{8q^4\cos^4(\theta/2)}{\pi^2(q-2)^2}\frac{(q^2-1)^2}{36} \\
        &= \frac{2q^4(q^2-1)^2\cos^4(\pi/q)}{9\pi^2(q-2)^2},
    \end{align*}
    as desired.
\end{proof}

\begin{cor}\label{cor:qveechpairlimit}
    We have
    \begin{equation*}
        \lim_{A \to \infty} \frac{c_{A}(\Vc_q)}{A} = 2\lambda_{\B_1^4} \cdot c_{\SV}(\Vc_q)^2.
    \end{equation*}
\end{cor}
\begin{proof}
    It is straightforward to verify that the single-variable Siegel-Veech constant $c_{\SV}(\Vc_q)$ is given by
    \begin{equation*}
        c_{\SV}(\Vc_q) = \frac{h}{\Covol(\Gamma)}\sum_{j=1}^p \frac{1}{\ell_j^2} = \frac{q^2(q^2-1)\cos^2(\pi/q)}{3\pi^2(q-2)},
    \end{equation*}
    and so simplifying yields the desired result.
\end{proof}

\begin{rmk}
    The asymptotic for Corollary \ref{cor:qveechpairlimit} can be established for a general Veech curve $\Vc$. This follows from the asymptotics of \cite[Section 7]{BurFai}.
\end{rmk}

\subsection{Branched covers of general Veech surfaces}

Let $\Vc^1 = \SL_2(\R)/\Gamma$ be a Veech curve based at $\omega \in \Hc^1$. For each $n \in \N$, there is a natural fiber bundle $\Ec_\Vc^n$ so that the fiber over a point $g\omega \in \Vc$ is $X^n$, where $X$ is the underlying flat surface of $\omega$, equipped with an affine $\Gamma$-action. We may alternatively describe
\begin{equation*}
    \Ec_\Vc^n = (\SL_2(\R) \times X^n)/\Gamma.
\end{equation*}
Note that $\Ec_{\Hc(0)}^n = \Hc^1(0,0^n)$. The space $\Ec_{\Vc}^n$ is a homogeneous space if and only if $\Vc = \Hc(0)$, and this homogeneity is used extensively in \cite{MarMoh}. While there are analogues of Ratner's theorems in the setting of general $\Ec_\Vc^n$ (see e.g. \cite{EskMarMor} and \cite{ChaiWeiYgo}), the analogue of Theorem \ref{thm:MarMoh} is not available in this setting. Moreover, for non-arithmetic $\Vc$ it is not clear what the set of saddle connections of even a \textit{generic} marked surface $\omega \in \Ec_\Vc^n$ should look like.

\begin{que}
    Suppose $\Vc$ is a Veech curve other than $\Hc(0)$. What is $c_A(\Ec_\Vc^n)$?
\end{que}

\subsection{Pair Siegel-Veech measures}\label{subsec:SVmeasveech}

We also record the pair Siegel-Veech measures associated to Veech curves, as defined in Section \ref{subsec:SVtransforms}, which can be readily derived from \cite[Theorem 1.8]{BurFai}.

\begin{prop}
    If $\Vc$ is a Veech curve, then $\eta_\Vc$ is some finitely supported measure on $\kay_\Vc$, and
    \begin{equation*}
        \nu_{\Vc}= \frac{1}{\Vol(\Vc^1)}\sum_{\substack{t \in \VT(\Vc) \\ t \ne 0}} \phi_\Vc(t)\bigl[\delta_{-t} + \delta_{+t}\bigr].
    \end{equation*}
\end{prop}

\printbibliography

@article{Fai, title={A higher moment formula for the Siegel–Veech transform over quotients by Hecke triangle groups}, volume={15}, ISSN={1661-7207, 1661-7215}, url={https://ems.press/doi/10.4171/ggd/591}, DOI={10.4171/ggd/591}, number={1}, journal={Groups, Geometry, and Dynamics}, author={Fairchild, Samantha}, year={2021}, pages={57--81} }

@article{MarMoh, title={Quantitative version of the Oppenheim conjecture for inhomogeneous quadratic forms}, volume={158}, ISSN={0012-7094, 1547-7398}, url={https://projecteuclid.org/journals/duke-mathematical-journal/volume-158/issue-1/Quantitative-version-of-the-Oppenheim-conjecture-for-inhomogeneous-quadratic-forms/10.1215/00127094-1276319.full}, DOI={10.1215/00127094-1276319}, number={1}, journal={Duke Mathematical Journal}, author={Margulis, Gregory and Mohammadi, Amir}, year={2011}, month=may, pages={121–160}, language={en} }

@article{EskMarMoz, title={Upper Bounds and Asymptotics in a Quantitative Version of the Oppenheim Conjecture}, volume={147}, ISSN={0003486X}, url={https://www.jstor.org/stable/120984?origin=crossref}, DOI={10.2307/120984}, number={1}, journal={The Annals of Mathematics}, author={Eskin, Alex and Margulis, Gregory and Mozes, Shahar}, year={1998}, month=jan, pages={93--141} }

@article{AthFaiMas, title={Counting pairs of saddle connections}, volume={431}, ISSN={00018708}, url={https://linkinghub.elsevier.com/retrieve/pii/S0001870823003766}, DOI={10.1016/j.aim.2023.109233}, journal={Advances in Mathematics}, author={Athreya, Jayadev and Fairchild, Samantha and Masur, Howard}, year={2023}, month=oct, pages={109233}, language={en} }

@article{AthCheMas, title={Siegel–Veech transforms are in $L^2$ (with an appendix by Jayadev S. Athreya and Rene Rühr)}, volume={14}, ISSN={1930-5311, 1930-532X}, url={https://www.aimsciences.org/article/doi/10.3934/jmd.2019001}, DOI={10.3934/jmd.2019001}, journal={Journal of Modern Dynamics}, author={Athreya, Jayadev S. and Cheung, Yitwah and Masur, Howard}, year={2019}, pages={1–19} }

@article{BurFai, title={Pairs in discrete lattice orbits with applications to Veech surfaces (with an appendix by Jon Chaika)}, ISSN={1435-9855, 1435-9863}, url={https://ems.press/doi/10.4171/jems/1563}, DOI={10.4171/jems/1563}, journal={Journal of the European Mathematical Society}, author={Burrin, Claire and Fairchild, Samantha}, year={2024}, month=nov }

@article{NevSte, title={Analogs of Wiener’s Ergodic Theorems for Semisimple Groups I}, volume={145}, ISSN={0003486X}, url={https://www.jstor.org/stable/2951845?origin=crossref}, DOI={10.2307/2951845}, number={3}, journal={The Annals of Mathematics}, author={Nevo, Amos and Stein, Elias M.}, year={1997}, month=may, pages={565--595} }

@article{ChaiWeiYgo, title={Horocycle dynamics in rank one invariant subvarieties. {I}: Weak measure classification and equidistribution}, journal={Israel Journal of Mathematics}, volume={271}, number={1}, pages={341--397}, DOI={10.1007/s11856-025-2805-3}, url={https://doi.org/10.1007/s11856-025-2805-3}, eprint={2301.12419}, eprinttype={arxiv}, eprintclass={math.DS}, author={Chaika, Jon and Weiss, Barak and Ygouf, Florent}, year={2026}, month=mar }

@article{EskMarMoz1, title={Quadratic forms of signature (2,2) and eigenvalue spacings on rectangular 2-tori}, volume={161}, ISSN={0003-486X}, url={http://annals.math.princeton.edu/2005/161-2/p03}, DOI={10.4007/annals.2005.161.679}, number={2}, journal={Annals of Mathematics}, author={Eskin, Alex and Margulis, Grigori and Mozes, Shahar}, year={2005}, month=mar, pages={679–725}, language={en} }

@article{Kim, title={Moments of Margulis functions and indefinite ternary quadratic forms}, url={http://arxiv.org/abs/2403.16563}, DOI={10.48550/arXiv.2403.16563}, note={arXiv:2403.16563 [math.DS]}, number={arXiv:2403.16563}, publisher={arXiv}, author={Kim, Wooyeon}, year={2025}, month=jul }

@article{BarShu, title={Unique ergodicity of branched covers of translation surfaces}, url={http://arxiv.org/abs/2606.05492}, DOI={10.48550/arXiv.2606.05492}, note={arXiv:2606.05492 [math.DS]}, number={arXiv:2606.05492}, publisher={arXiv}, author={Baron, Polina and Shuvaeva, Elizaveta}, year={2026}, month=jun}

@article{EskMasSch, title={Billiards in rectangles with barriers}, volume={118}, ISSN={0012-7094}, url={https://projecteuclid.org/journals/duke-mathematical-journal/volume-118/issue-3/Billiards-in-rectangles-with-barriers/10.1215/S0012-7094-03-11832-3.full}, DOI={10.1215/S0012-7094-03-11832-3}, number={3}, journal={Duke Mathematical Journal}, author={Eskin, Alex and Masur, Howard and Schmoll, Martin}, year={2003}, month=jun, pages={427--463}}

@article{EskMarMor, title={Unipotent flows on the space of branched covers of Veech surfaces}, volume={26}, rights={https://www.cambridge.org/core/terms}, ISSN={0143-3857, 1469-4417}, url={https://www.cambridge.org/core/product/identifier/S0143385705000234/type/journal_article}, DOI={10.1017/S0143385705000234}, number={1}, journal={Ergodic Theory and Dynamical Systems}, author={Eskin, Alex and Marklof, Jens and Morris, Dave Witte}, year={2006}, month=feb, pages={129–162}, language={en} }

@article{ElBMarVin, title={The Distribution of Directions in an Affine Lattice: Two-Point Correlations and Mixed Moments}, ISSN={1073-7928, 1687-0247}, url={https://academic.oup.com/imrn/article-lookup/doi/10.1093/imrn/rnt258}, DOI={10.1093/imrn/rnt258}, journal={International Mathematics Research Notices}, author={El-Baz, D. and Marklof, J. and Vinogradov, I.}, year={2015}, volume={2015}, number={5}, month=may, pages={1371--1400}, language={en} }

@book{Ten, address={Providence, Rhode Island}, series={Graduate Studies in Mathematics}, number={163}, edition={3}, title={Introduction to analytic and probabilistic number theory}, ISBN={9781470422233}, publisher={American Mathematical Society}, author={Tenenbaum, Gérald}, translator={Ion, Patrick D. F.}, year={2015}, collection={Graduate Studies in Mathematics}, language={eng} }

@book{Sza, address={Cambridge}, series={Cambridge studies in advanced mathematics}, title={Galois groups and fundamental groups}, ISBN={9780521888509}, callNumber={512.32}, publisher={Cambridge University Press}, author={Szamuely, Tamás}, year={2009}, collection={Cambridge studies in advanced mathematics}, language={eng} }

@article{Ore, title={Some remarks on commutators}, volume={2}, ISSN={1088-6826, 0002-9939}, url={https://www.ams.org/proc/1951-002-02/S0002-9939-1951-0040298-4/}, DOI={10.1090/S0002-9939-1951-0040298-4}, number={2}, journal={Proceedings of the American Mathematical Society}, author={Ore, Oystein}, year={1951}, month=apr, pages={307–314}, language={en} }

@unpublished{FaiFenKatMas,
    label = {FFKM},
  author = {Fairchild, Samantha and Fenves, Inbo and Katipearachchi, Kanishka and Masur, Howard},
  title  = {Pair Siegel-Veech constants for strata},
  note   = {In preparation},
}

@article{EskMasZor, title={Moduli spaces of abelian differentials: the principal boundary, counting problems, and the Siegel-Veech constants}, volume={97}, rights={http://www.springer.com/tdm}, ISSN={0073-8301, 1618-1913}, url={https://pmihes.centre-mersenne.org/articles/10.1007/s10240-003-0015-1/}, DOI={10.1007/s10240-003-0015-1}, journal={Publications Mathématiques de l’IHÉS}, author={Eskin, Alex and Masur, Howard and Zorich, Anton}, year={2003}, month=sep, pages={61–179}, language={en} }

@article{CheMolSauZag, label={CMSZ}, title={Masur–Veech volumes and intersection theory on moduli spaces of Abelian differentials}, volume={222}, ISSN={0020-9910, 1432-1297}, url={https://link.springer.com/10.1007/s00222-020-00969-4}, DOI={10.1007/s00222-020-00969-4}, number={1}, journal={Inventiones mathematicae}, author={Chen, Dawei and Möller, Martin and Sauvaget, Adrien and Zagier, Don}, year={2020}, month=oct, pages={283–373}, language={en} }

@article{Sch, title={On the asymptotic quadratic growth rate of saddle connections and periodic orbits on marked flat tori}, volume={12}, rights={http://www.springer.com/tdm}, ISSN={1016-443X, 1420-8970}, url={http://link.springer.com/10.1007/s00039-002-8260-x}, DOI={10.1007/s00039-002-8260-x}, number={3}, journal={Geometric And Functional Analysis}, author={Schmoll, M.}, year={2002}, month=aug, pages={622–649} }

@article{Sch1, title={Moduli spaces of branched covers of Veech surfaces I: d-symmetric differentials}, url={http://arxiv.org/abs/math/0602396}, DOI={10.48550/arXiv.math/0602396}, note={arXiv:math/0602396}, number={arXiv:math/0602396}, publisher={arXiv}, author={Schmoll, Martin}, year={2006}, month=feb }

@article{AulCalMatSalSch, label={ACMSS}, title={Siegel–Veech constants for cyclic covers of generic translation surfaces}, journal={Journal of the London Mathematical Society}, volume={113}, number={1}, eid={e70413}, DOI={10.1112/jlms.70413}, url={https://doi.org/10.1112/jlms.70413}, eprint={2409.06600}, eprinttype={arxiv}, eprintclass={math.DS}, author={Aulicino, David and Calderon, Aaron and Matheus, Carlos and Salter, Nick and Schmoll, Martin}, year={2026}, month=jan }

@article{Vee, title={Teichm\"uller curves in moduli space, Eisenstein series and an application to triangular billiards}, volume={97}, rights={http://www.springer.com/tdm}, ISSN={0020-9910, 1432-1297}, url={http://link.springer.com/10.1007/BF01388890}, DOI={10.1007/BF01388890}, number={3}, journal={Inventiones Mathematicae}, author={Veech, W. A.}, year={1989}, month=oct, pages={553–583}, language={de} }

@article{Ngu, title={Intersection theory and Siegel-Veech constants for Prym eigenform loci in $\Omega\mathcal{M}_3(2,2)^{odd}$}, url={http://arxiv.org/abs/2510.23333}, DOI={10.48550/arXiv.2510.23333}, note={arXiv:2510.23333 [math.GT]}, number={arXiv:2510.23333}, publisher={arXiv}, author={Nguyen, Duc-Manh}, year={2025}, month=dec }

@article{Sie, title={A Mean Value Theorem in Geometry of Numbers}, volume={46}, ISSN={0003486X}, url={https://www.jstor.org/stable/1969027?origin=crossref}, DOI={10.2307/1969027}, number={2}, journal={The Annals of Mathematics}, author={Siegel, Carl Ludwig}, year={1945}, month=apr, pages={340--347} }

@article{Gou, title={Siegel–Veech Constants for Strata of Moduli Spaces of Quadratic Differentials}, volume={25}, ISSN={1016-443X, 1420-8970}, url={http://link.springer.com/10.1007/s00039-015-0345-4}, DOI={10.1007/s00039-015-0345-4}, number={5}, journal={Geometric and Functional Analysis}, author={Goujard, Elise}, year={2015}, month=oct, pages={1440–1492}, language={en} }

@article{DelGouZogZor, label={DGZZ}, title={Masur–Veech volumes, frequencies of simple closed geodesics, and intersection numbers of moduli spaces of curves}, volume={170}, ISSN={0012-7094}, url={https://projecteuclid.org/journals/duke-mathematical-journal/volume-170/issue-12/MasurVeech-volumes-frequencies-of-simple-closed-geodesics-and-intersection-numbers/10.1215/00127094-2021-0054.full}, DOI={10.1215/00127094-2021-0054}, number={12}, journal={Duke Mathematical Journal}, author={Delecroix, Vincent and Goujard, Élise and Zograf, Peter and Zorich, Anton}, year={2021}, month=sep, pages={2633--2718} }

@article{AnaMulNor, title={Weil-Petersson volumes, stability conditions and wall-crossing}, url={http://arxiv.org/abs/2310.13281}, DOI={10.48550/arXiv.2310.13281}, note={arXiv:2310.13281 [math.AG]}, number={arXiv:2310.13281}, publisher={arXiv}, author={Anagnostou, Lukas and Mullane, Scott and Norbury, Paul}, year={2023}, month=oct }

@article{GutHubSch, title={Affine diffeomorphisms of translation surfaces: Periodic points, Fuchsian groups, and arithmeticity}, volume={36}, rights={http://www.elsevier.com/tdm/userlicense/1.0/}, ISSN={00129593}, url={http://linkinghub.elsevier.com/retrieve/pii/S0012959303000491}, DOI={10.1016/j.ansens.2003.05.001}, number={6}, journal={Annales Scientifiques de l’École Normale Supérieure}, author={Gutkin, E and Hubert, P and Schmidt, T}, year={2003}, month=nov, pages={847–866}, language={en} }

@article{Schm, title={A metrical theorem in geometry of numbers}, volume={95}, ISSN={1088-6850, 0002-9947}, url={https://www.ams.org/tran/1960-095-03/S0002-9947-1960-0117222-9/}, DOI={10.1090/S0002-9947-1960-0117222-9}, number={3}, journal={Transactions of the American Mathematical Society}, author={Schmidt, Wolfgang M.}, year={1960}, pages={516–529}, language={en} }

@article{EskMcM, title={Mixing, counting, and equidistribution in Lie groups}, volume={71}, ISSN={0012-7094, 1547-7398}, url={https://projecteuclid.org/journals/duke-mathematical-journal/volume-71/issue-1/Mixing-counting-and-equidistribution-in-Lie-groups/10.1215/S0012-7094-93-07108-6.full}, DOI={10.1215/S0012-7094-93-07108-6}, number={1}, journal={Duke Mathematical Journal}, author={Eskin, Alex and McMullen, Curt}, year={1993}, month=jul, pages={181–209} }

@article{SmiWei, title={Characterizations of lattice surfaces}, volume={180}, rights={http://www.springer.com/tdm}, ISSN={0020-9910, 1432-1297}, url={http://link.springer.com/10.1007/s00222-010-0236-0}, DOI={10.1007/s00222-010-0236-0}, number={3}, journal={Inventiones mathematicae}, author={Smillie, John and Weiss, Barak}, year={2010}, month=jun, pages={535–557}, language={en} }

@article{Fil, title={Relative periods of abelian differentials in a prescribed stratum}, url={http://arxiv.org/abs/2503.21285}, DOI={10.48550/arXiv.2503.21285}, note={arXiv:2503.21285 [math.GT]}, number={arXiv:2503.21285}, publisher={arXiv}, author={Le Fils, Thomas}, year={2025}, month=mar }

@article{ArnHub, title={Fractions continues sur les surfaces de Veech}, volume={81}, rights={http://www.springer.com/tdm}, ISSN={0021-7670, 1565-8538}, url={http://link.springer.com/10.1007/BF02788985}, DOI={10.1007/BF02788985}, number={1}, journal={Journal d’Analyse Mathématique}, author={Arnoux, Pierre and Hubert, Pascal}, year={2000}, month=dec, pages={35–64}, language={fr} }

@article{ArnSch, title={Veech surfaces with nonperiodic directions in the trace field}, volume={3}, rights={http://creativecommons.org/licenses/by/3.0/}, ISSN={1930-5311}, url={https://www.aimsciences.org/en/article/doi/10.3934/jmd.2009.3.611}, DOI={10.3934/jmd.2009.3.611}, number={4}, journal={Journal of Modern Dynamics}, publisher={Journal of Modern Dynamics}, author={Arnoux, Pierre and Schmidt, Thomas A.}, year={2009}, pages={611–629}, language={en} }

@book{GorNev, address={Princeton, N.J}, series={Annals of mathematics studies}, title={The ergodic theory of lattice subgroups}, ISBN={9780691141848}, callNumber={QA313 .G67 2010}, publisher={Princeton University Press}, author={Gorodnik, Alexander and Nevo, Amos}, year={2010}, collection={Annals of mathematics studies} }

@article{Vee1, title={Siegel Measures}, volume={148}, ISSN={0003486X}, url={https://www.jstor.org/stable/121033?origin=crossref}, DOI={10.2307/121033}, number={3}, journal={The Annals of Mathematics}, author={Veech, William A.}, year={1998}, month=nov, pages={895--944} }

@article{SzaTur,
  author  = {Szalay, M. and Tur{\'a}n, P.},
  title   = {On some problems of the statistical theory of partitions
             with application to characters of the symmetric group. {I}},
  journal = {Acta Mathematica Academiae Scientiarum Hungaricae},
  volume  = {29},
  number  = {3--4},
  pages   = {361--379},
  year    = {1977},
  doi     = {10.1007/BF01895857}
}

@article{BriFraGar, title={Asymptotics for the twisted eta-product and applications to sign changes in partitions}, volume={9}, ISSN={2197-9847}, url={https://doi.org/10.1007/s40687-022-00355-x}, DOI={10.1007/s40687-022-00355-x}, number={4}, journal={Research in the Mathematical Sciences}, author={Bridges, Walter and Franke, Johann and Garnowski, Taylor}, year={2022}, month=sep, pages={61}, language={en} }

@book{Goo, address={Berlin, Heidelberg}, series={Lecture Notes in Mathematics}, title={Local Analysis of Selberg’s Trace Formula}, volume={1040}, rights={http://www.springer.com/tdm}, ISBN={9783540127130}, url={http://link.springer.com/10.1007/BFb0073074}, DOI={10.1007/BFb0073074}, publisher={Springer}, author={Good, Anton}, year={1983}, collection={Lecture Notes in Mathematics} }

@misc{StrSodVis,
  author        = {Str{\"o}mbergsson, Andreas and S{\"o}dergren, Anders
                   and Vishe, Pankaj},
  title         = {Effective equidistribution of unipotent orbits in
                   homogeneous spaces of
                   {$\mathrm{SL}(2,\mathbb{R})\ltimes(\mathbb{R}^{2})^{k}$}},
  year          = {2026},
  eprint        = {2604.08753},
  archivePrefix = {arXiv},
  primaryClass  = {math.DS},
  doi           = {10.48550/arXiv.2604.08753},
  url           = {https://arxiv.org/abs/2604.08753}
}

@incollection{Gut,
  author    = {Gutkin, Eugene},
  title     = {Branched coverings and closed geodesics in flat surfaces,
               with applications to billiards},
  booktitle = {Dynamical Systems: From Crystal to Chaos},
  editor    = {Gambaudo, J.-M. and Hubert, P. and
               Tisseur, P. and Vaienti, S.},
  publisher = {World Scientific},
  address   = {Singapore},
  year      = {2000},
  pages     = {259--273},
  doi       = {10.1142/9789812793829_0024},
  url       = {https://doi.org/10.1142/9789812793829_0024}
}

@article{HubSch,
  author  = {Hubert, Pascal and Schmidt, Thomas A.},
  title   = {Invariants of translation surfaces},
  journal = {Annales de l'Institut Fourier},
  volume  = {51},
  number  = {2},
  year    = {2001},
  pages   = {461--495},
  doi     = {10.5802/aif.1829},
  url     = {https://numdam.org/articles/10.5802/aif.1829/}
}

@article{San, title={Gaps of saddle connection directions for some branched covers of tori}, ISSN={0143-3857, 1469-4417}, url={https://www.cambridge.org/core/product/identifier/S014338572100081X/type/journal_article}, DOI={10.1017/etds.2021.81}, journal={Ergodic Theory and Dynamical Systems}, author={Sanchez, Anthony}, year={2022}, volume={42}, number={10}, month=oct, pages={3191--3245}, language={en} }

@book{AthMas, address={Providence, Rhode Island}, edition={First edition}, series={Graduate studies in mathematics}, title={Translation surfaces}, ISBN={9781470476779}, callNumber={QA333 .A84 2024}, publisher={American Mathematical Society}, author={Athreya, Jayadev S. and Masur, Howard}, year={2024}, collection={Graduate studies in mathematics} }

@article{EskMas,
  title = {Asymptotic formulas on flat surfaces},
  volume = {21},
  rights = {https://www.cambridge.org/core/terms},
  issn = {0143-3857, 1469-4417},
  url = {https://www.cambridge.org/core/product/identifier/S0143385701001225/type/journal_article},
  doi = {10.1017/S0143385701001225},
  number = {2},
  journal = {Ergodic Theory and Dynamical Systems},
  author = {Eskin, Alex and Masur, Howard},
  year = {2001},
  month = apr,
  pages = {443–478},
  language = {en},
}

@article{McM, title={Teichmüller geodesics of infinite complexity}, volume={191}, ISSN={0001-5962}, url={http://projecteuclid.org/euclid.acta/1485891606}, DOI={10.1007/BF02392964}, number={2}, journal={Acta Mathematica}, author={McMullen, Curtis T.}, year={2003}, pages={191–223}, language={en} }

@article{Mar, title={Pair correlation densities of inhomogeneous quadratic forms}, volume={158}, ISSN={0003-486X}, url={http://annals.math.princeton.edu/2003/158-2/p02}, DOI={10.4007/annals.2003.158.419}, number={2}, journal={Annals of Mathematics}, author={Marklof, Jens}, year={2003}, month=sep, pages={419–471}, language={en} }

@article{KenSmi, title={Billiards on rational-angled triangles}, volume={75}, ISSN={0010-2571, 1420-8946}, url={https://ems.press/doi/10.1007/s000140050113}, DOI={10.1007/s000140050113}, number={1}, journal={Commentarii Mathematici Helvetici}, author={Kenyon, Richard and Smillie, John}, year={2000}, month=mar, pages={65–108} }

\end{document}